%% file: 2003-12.tex
\documentclass{gtart}
\input gtoutput

\lognumber{281}
\volumenumber{7}\papernumber{12}
\volumeyear{2003}
\pagenumbers{399}{441}
\received{30 September 2002}
\accepted{13 April 2003}
\published{24 June 2003}
\proposed{Jean-Pierre Otal}
\seconded{Walter Neumann, Martin Bridson}

\usepackage{amsmath,amssymb}
\usepackage{graphicx,psfrag}

\newcommand{\C}{{\mathbb C}}
\newcommand{\Q}{{\mathbb Q}}
\newcommand{\F}{{\mathbb F}}
\newcommand{\R}{{\mathbb R}}
\newcommand{\Z}{{\mathbb Z}}
\newcommand{\T}{{\mathcal T}}
\renewcommand{\H}{{\mathbb H}}
\renewcommand{\O}{{\mathcal{O}}}
\newcommand{\p}{\wp}
\newcommand{\bdry}{\partial}
\newcommand{\iso}{\cong}

\newcommand{\maps}{\colon\thinspace}

\DeclareMathOperator{\tr}{tr}

\DeclareMathOperator{\Isom}{Isom}
\DeclareMathOperator{\Teich}{Teich}

\newcommand{\CH}{{ \mathbb{C}\, \mathrm{H}}}
\newcommand{\RP}{{ \mathbb{R}\, \mathrm{P}}}
\newcommand{\PSL}[2]{\mathrm{PSL}_{#1} #2}
\newcommand{\PGL}[2]{\mathrm{PGL}_{#1} #2}
\newcommand{\SL}[2]{\mathrm{SL}_{#1} #2}
\newcommand{\GL}[2]{\mathrm{GL}_{#1} #2}
\newcommand{\abs}[1]{{\left| #1 \right|}}
\newcommand{\pair}[1]{\left\langle #1 \right\rangle}
\newcommand{\unfilled}{\, \cdot \, }
\newcommand{\spandef}[2]{{  \left\langle  {#1}  \ \left| \   {#2} \right. \right\rangle }}
\newcommand{\setdef}[2]{{  \left\{  {#1}  \ \left| \   {#2} \right. \right\} }}
\newcommand{\mtext}[1]{\quad\mbox{#1}\quad}

\newcommand{\SnapPea}{\texttt{SnapPea}}
\newcommand{\GAP}{\texttt{GAP}}
\newcommand{\Snap}{\texttt{Snap}}
\newcommand{\MAGMA}{\texttt{MAGMA}}

\theoremstyle{plain} 

\swapnumbers
\newtheorem{theorem}{\ Theorem}[section]
\newtheorem{conjecture}[theorem]{\ Conjecture}
\newtheorem{lemma}[theorem]{\ Lemma}
\newtheorem{corollary}[theorem]{\ Corollary}

\newtheorem{claim}[theorem]{\ Claim}
\newtheorem{question}[theorem]{\ Question}

\newtheorem{virtual_haken}[theorem]{\ Virtual Haken Conjecture}
\newtheorem{hyperbolic_dehn_surgery}[theorem]{\ Hyperbolic Dehn Surgery Theorem}
\newtheorem{vpbn}[theorem]{\ Virtual Positive Betti Number
  Conjecture}

\theoremstyle{definition}

\theoremstyle{remark}

\swapnumbers
\theoremstyle{plain}\newtheorem*{whiteforward}{Theorem (\ref{whiteheadthm})}

\newtheoremstyle{bold}{14pt plus6pt minus6pt}%
{2pt plus3pt minus3pt}{\bf}{}{\bf}{}{.66em}%
{\thmname{#1}\thmnumber{#2}\thmnote{\bf\stdspace[#3]}}
\theoremstyle{bold}
\newtheorem{sshead}[theorem]{}

\def\subsection#1{\penalty -800\begin{sshead}#1\end{sshead}\penalty 800}
\def\sh#1{\par\medskip{\bf #1}\par\vskip 2pt plus3pt minus3pt}

\begin{document}

\title{The virtual Haken conjecture:\\Experiments and examples}

\author{Nathan M Dunfield\\William P Thurston}
\shortauthors{Dunfield and Thurston}

\address{Department of Mathematics,  Harvard University\\Cambridge MA, 02138, USA} 
\email{nathand@math.harvard.edu}

\secondaddress{Department of Mathematics, University of California, Davis \\  Davis, CA 95616, USA}
\secondemail{wpt@math.ucdavis.edu}

\asciiaddress{Department of Mathematics,  Harvard University\\Cambridge MA, 02138, USA\\and\\Department of Mathematics, University of California, 
Davis\\Davis, CA 95616, USA}

\asciiemail{nathand@math.harvard.edu, wpt@math.ucdavis.edu}


\begin{abstract} 
  A 3-manifold is Haken if it contains a topologically essential
  surface.  The Virtual Haken Conjecture says that every irreducible
  3-manifold with infinite fundamental group has a finite cover which
  is Haken.  Here, we discuss two interrelated topics concerning
  this conjecture.
  
  First, we describe computer experiments which give strong evidence
  that the Virtual Haken Conjecture is true for hyperbolic
  3-manifolds. We took the complete Hodgson-Weeks census of 10,986
  small-volume closed hyperbolic 3-manifolds, and for each of them
  found finite covers which are Haken.  There are interesting and
  unexplained patterns in the data which may lead to a better
  understanding of this problem.
  
  Second, we discuss a method for transferring the virtual Haken
  property under Dehn filling.  In particular, we show that if a
  3-manifold with torus boundary has a Seifert fibered Dehn filling
  with hyperbolic base orbifold, then most of the Dehn filled
  manifolds are virtually Haken.  We use this to show that every
  non-trivial Dehn surgery on the figure-8 knot is virtually Haken.
\end{abstract}

{\small\maketitlepage}

\section{Introduction}

Let $M$ be an orientable 3-manifold.  A properly embedded orientable
surface $S \neq S^2$ in $M$ is \emph{incompressible} if it is not
boundary parallel, and the inclusion $\pi_1(S) \to \pi_1(M)$ is
injective.  A manifold is \emph{Haken} if it is irreducible and
contains an incompressible surface.  Haken manifolds are by far the
best understood class of 3-man\-i\-folds.  This is because splitting a
Haken manifold along an incompressible surface results in a simpler
Haken manifold.  This allows induction arguments for these manifolds.

However, many irreducible 3-manifolds with infinite fundamental group
are not Haken, e.g.~all but 4 Dehn surgeries on the
figure-8 knot.  It has been very hard to prove anything about
non-Haken manifolds, at least without assuming some sort of additional
Haken-like structure, such as a foliation or lamination. 

Sometimes, a non-Haken 3-manifold $M$ has a finite cover which is Haken.
Most of the known properties for Haken manifolds can then be pushed
down to $M$ (though showing this can be difficult).  Thus, one of the
most interesting conjectures about 3-manifolds is Waldhausen's
conjecture \cite{Waldhausen68}:
\begin{virtual_haken}
  Suppose $M$ is an irreducible 3-manifold with infinite fundamental
  group.  Then $M$ has a finite cover which is Haken.  
\end{virtual_haken}
A 3-manifold satisfying this conjecture is called \emph{virtually
  Haken}.  For more background and references on this conjecture see
Kirby's problem list \cite{KirbyList}, problems 3.2, 3.50, and 3.51.
See also \cite{CooperLong99, CooperLongSSSSS} and \cite{Lubotzky96,
  Lubotzky98} for some of the latest results toward this conjecture.
The importance of this conjecture is enhanced because it's now known
that 3-manifolds which are virtually Haken are geometrizable
\cite{GabaiMeyerhoffThurston, Gabai97,Scott83,
  MeeksSimonYau,Mess,Gabai92,CassonJungreis}.

There are several stronger forms of this conjecture, including asking
that the finite cover be not just Haken but a surface bundle over the
circle.  We will be interested in the following version.  Let $M$ be a
closed irreducible 3-manifold.  If $H_2(M,\Z) \neq 0$ then $M$ is
Haken, as any non-zero class in $H_2(M,\Z)$ can be represented by an
incompressible surface.  Now $H_2(M, \Z)$ is isomorphic to $H^1(M,\Z)$
by Poincar{\'e} duality, and $H^1(M,\Z)$ is a free abelian group.  So if
the first betti-number of $M$ is $\beta_1(M) = \dim H_1(M, \R) =\dim
H^1(M, \R)$, then $\beta_1(M) > 0$ implies $M$ is Haken. As the cover
of an irreducible 3-manifold is irreducible \cite{MeeksSimonYau}, a
stronger form of the Virtual Haken Conjecture is:
\begin{vpbn}\label{vpbn}
  Suppose $M$ is an irreducible 3-mani\-fold with infinite fundamental
  group.  Then $M$ has a finite cover $N$ where $\beta_1(N) > 0$.
\end{vpbn}
We will say that such an $M$ has \emph{virtual positive betti number}.
Note that $\beta_1(N) > 0$ if and only if $H_1(N, \Z)$, the
abelianization of $\pi_1(N)$, is infinite.  So an equivalent, more
algebraic, formulation of Conjecture~\ref{vpbn} is:
\begin{conjecture}
  Suppose $M$ is an irreducible 3-manifold.  Assume that $\pi_1(M)$ is
  infinite.  Then $\pi_1(M)$ has a finite index subgroup with infinite
  abelianization.
\end{conjecture}
Here, we focus on this form of the Virtual Haken Conjecture because
its algebraic nature makes it easier to examine both theoretically and
computationally.  While in theory one can to use normal surface
algorithms to decide if a manifold is Haken \cite{JacoOertel}, in
practice these algorithms are prohibitively slow in all but the
simplest examples.  Computing homology is much easier as it boils down
to computing the rank of a matrix.  Also, it's probably true that
having virtual positive betti number isn't much stronger than being
virtually Haken (see the discussion of \cite{Lubotzky96} in
Section~\ref{dehn_questions} below).

\subsection{Outline of the paper}

This paper examines the Virtual Haken Conjecture in two interrelated
parts:

\sh{Experiment: Sections \ref{experiment}-\ref{questions}}   

Here, we describe experiments which strongly support the Virtual
Positive Betti Number Conjecture.  We looked at the 10,986
small-volume hyperbolic manifolds in the Hodgson-Weeks census, and
tried to show that they had virtual positive betti number.  In all
cases, we succeeded.  It was natural to restrict to hyperbolic
3-manifolds for our experiment since, in practice, all 3-manifolds
are geometrizable and the Virtual Positive Betti Number Conjecture is
known for all other kinds of geometrizable 3-manifolds.

Section~\ref{experiment} gives an overview of the experiment and
discusses the results and limitations of the survey.
Sections~\ref{techniques} and \ref{comp_rank} describe the techniques
used to compute the homology of the covers.  Section~\ref{simple}
discusses some interesting patterns that we found among the covers
where the covering group is a simple group.  Some further questions are 
given in Section~\ref{questions}.

\sh{Examples and Dehn filling: Sections \ref{transfer} - \ref{sister}}
Here we consider Dehn fillings of a fixed 3-man\-i\-fold $M$ with torus
boundary.  Generalizing work of Boyer and Zhang \cite{BoyerZhang00},
we give a method for transferring virtual positive betti number from
one filling of $M$ to another.  Roughly, Theorem~\ref{SF_filling} says
that if $M$ has a filling which is Seifert fibered with hyperbolic
base orbifold, then most Dehn fillings have virtual positive betti
number.  We use this to give new examples of manifolds $M$ where all
but finitely many Dehn fillings have virtual positive betti number.
In Section~\ref{Whitehead}, we show this holds for most surgeries on
one component of the Whitehead link.

In the case of figure-8 knot, we use work of Holt and Plesken
\cite{HoltPlesken92} to amplify our results, and prove that every
non-trivial Dehn surgery on the figure-8 knot has virtual positive
betti number (Theorem~\ref{fig8thm}).

In Section~\ref{dehn_questions}, we discuss possible avenues to
other results using fillings which are Haken rather than Seifert
fibered.  This approach is easiest in the case of toroidal Dehn
fillings, and using these techniques we prove
(Theorem~\ref{sister_thm}) that all Dehn fillings on the sister of
the figure-8 complement satisfy the Virtual Positive Betti Number
Conjecture.

\sh{Acknowledgments} 
The first author was partially supported by an NSF Postdoctoral
Fellowship.  The second author was partially supported by NSF grants
DMS-9704135 and DMS-0072540. We would like to thank Ian Agol, Daniel
Allcock, Matt Baker, Danny Calegari, Greg Kuperberg, Darren Long, Alex
Lubotzky, Alan Reid, William Stein, and Dylan Thurston for useful
conversations.  We also thank of the authors of the computer programs
\SnapPea\ \cite{SnapPea} and \GAP\ \cite{GAP4} which were critical for
our computations.

\section{The experiment}\label{experiment}

\subsection{The manifolds}  
We looked at the 10,986 hyperbolic 3-manifolds in the Hodgson-Weeks
census of small-volume closed hyperbolic 3-manifolds \cite{SnapPea}.
The volumes of these manifolds range from that of the smallest known
manifold ($0.942707...$) to $6.5$.  While there are infinitely
many closed hyperbolic 3-manifolds with volume less than $6.5$, there
are only finitely many if we also bound the injectivity radius from
below.  The census manifolds are an approximation to all closed
hyperbolic 3-manifolds with volume $<6.5$ and injectivity radius $>
0.3$.

A more precise description of these manifolds is this.  Start with the
Callahan-Hildebrand-Weeks census of cusped finite-volume hyperbolic
3-manifolds, which is a complete list of the those having ideal
triangulations with 7 or fewer tetrahedra \cite{HildebrandWeeks,
  CallahanHildebrandWeeks}.  The closed census consists of all the
Dehn fillings on the 1-cusped manifolds in the cusped census, where
the closed manifold has shortest geodesic of length $>0.3$.

Only 132 of the 10,986 manifolds have positive betti number.  It is
also worth mentioning that many (probably the vast majority) of these
manifolds are non-Haken. For the 246 manifolds with volume less
than 3, exactly 15 are Haken \cite{Dunfield:haken}.

\subsection{Computational framework}

For each 3-manifold,  we started with a finite presentation of its
fundamental group $G$, and then looked for a finite index subgroup $H$
of $G$ which has infinite abelianization.  There is a fair amount of
literature on how find such an $H$, because finding a finite index
subgroup with infinite abelianization is one of the main computational
techniques for proving that a given finitely presented group is
infinite.  See \cite{Plesken99} for a survey.  The key idea which
simplifies the computations is contained in \cite{HoltPlesken92}, which
we used in the form described in Section~\ref{techniques}.

We used \SnapPea\ \cite{SnapPea} to give presentations for the
fundamental groups of each of the manifolds in the closed census.  We
then used \GAP\ \cite{GAP4} to find various finite index subgroups and
compute the homology of the subgroups (see
Sections~\ref{techniques}-\ref{comp_rank}).

\subsection{Types of covers}

When looking for a subgroup with positive betti number, we tried a
number of different types of subgroups.  Some types were much better
at producing homology than others.  Those that worked well were:
\begin{itemize}
\item Abelian/$p$-group covers with exponent $2$ or $3$.
  
\item Low ($<20$) index subgroups.  Coset enumeration techniques allow
  one to enumerate low-index subgroups \cite{Sims94}.  Given such a
  subgroup $H < G$, we looked at the largest normal subgroup
  contained in $H$, to maximize the chance of finding homology.
  
\item Normal subgroups where the quotient is a finite simple group.
  These were found by choosing the simple group in advance and then
  finding all epimorphisms of $G$ onto that group.

\end{itemize}
The following types were inefficient in producing homology:
\begin{itemize}
\item Abelian/nilpotent covers with exponents $> 3$.
\item Dihedral covers.
\item Intersections of subgroups of the types listed in the first list
  (the useful types).
\end{itemize}
It would be nice to have heuristics which explain why some things
worked and others didn't (we plan to explore this further in
\cite{DunfieldThurstonHeuristics}).  Also, while intersecting
subgroups was not efficient in general, there were certain manifolds
where the only positive betti number cover we could find were of this
type.

\subsection{Results}

We were able to find positive betti number covers for all of the
Hodgson-Weeks census manifolds.  For most of the manifolds, it was
easy to find such a cover.  For instance, just looking at abelian
covers and subgroups of index $\leq 6$ works for $42\%$ of the
manifolds.  See Table~\ref{quick_success} for more about the degrees
of the covers we used.
\begin{table}[tbhp]
\centering
\hspace{0.12\textwidth}\parbox[b]{0.18\textwidth}{
\begin{tabular}{rr}
  \textbf{$d$} & \textbf{\%}  \\
  \hline
  1      &  1.2  \\
  2      &  3.8  \\
  5      & 21.2  \\ 
  6      & 39.3  \\
  10     & 57.9  \\
  20     & 68.3  \\ 
  50     & 88.8  \\
 100     & 95.6  \\
 200     & 98.1  
\end{tabular}
}\hspace{0.1\textwidth}%
\parbox[c]{0.6\textwidth}{
\psfrag{x}{$\log(d)$}
\psfrag{y}{$\%$}
\includegraphics[scale=0.65]{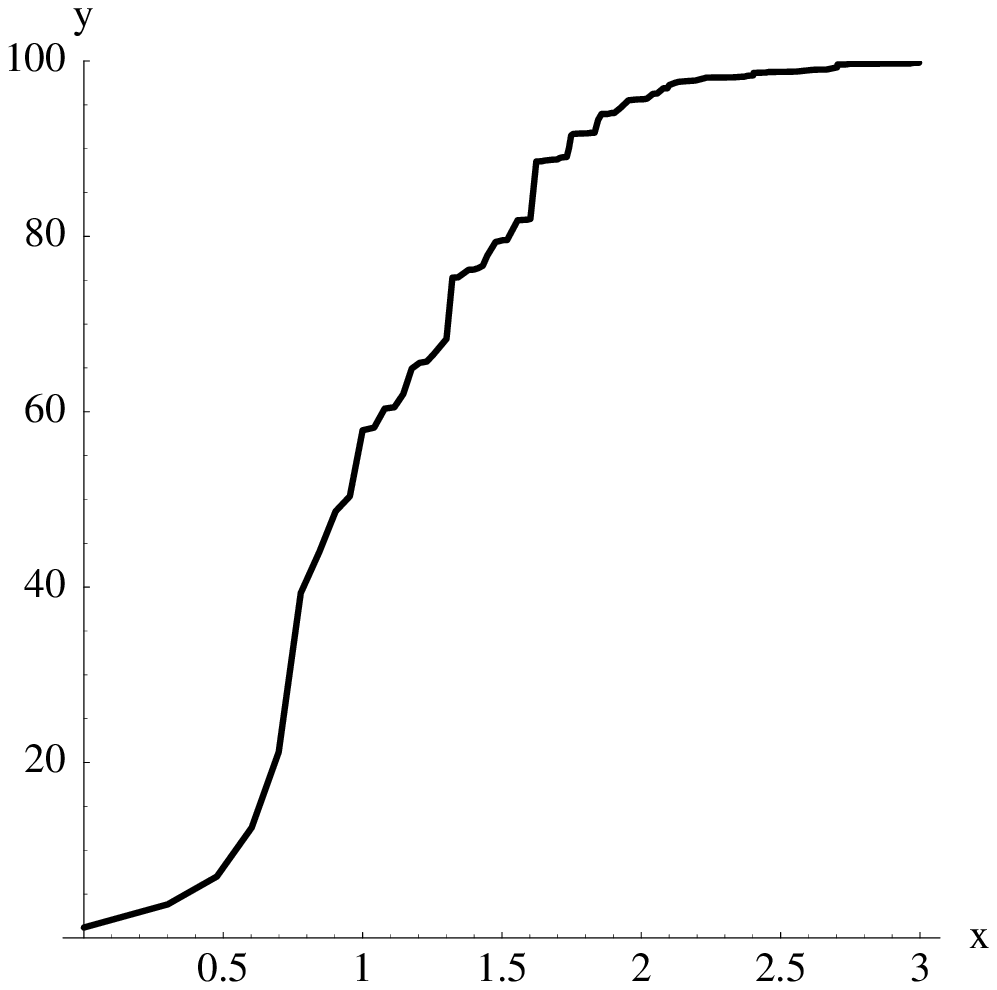}
}

\caption{
The table at left shows the proportion of manifolds for which we found a cover
with positive betti number of degree $\leq d$. Note this is just for the
covers that we found, which are not always the positive betti number
covers of smallest degree.   The plot at right presents all of the data,
where $\log(d)$ is base 10.
}
\label{quick_success}

\end{table}

For each of the manifolds, we stored a presentation of the fundamental
group and a homomorphism from that finitely presented group to $S_n$
whose kernel has positive betti number.  This information is available
on the web at \cite{VirtHakenWebsite} together with the \GAP\ code we
used for the computations, and will hopefully be useful as a source of
examples.  The amount of computer time used to find all the covers was
in excess of one CPU-year, but the amount of time needed to check all
the covers for homology given the data available at
\cite{VirtHakenWebsite} is only a few of hours.

\begin{figure}[ht!]
\begin{center}
\includegraphics[scale=0.3]{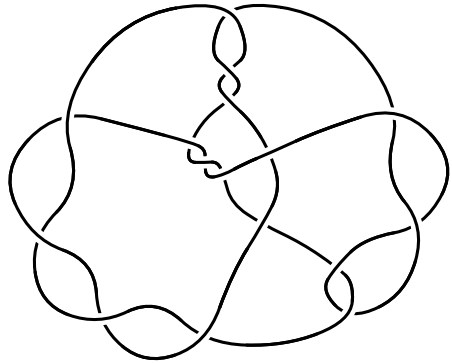}
\end{center}
\vspace{-0.3in}
\caption{
   The 2-fold cover branched over this knot is the manifold $N$.
   Figure created with \cite{Knotscape}.}
\label{knot}
\end{figure}

There was one manifold in particular where it was very difficult to
find a cover with positive betti number.  This manifold is $N =
s633(2,3)$.  Its volume is $4.49769817315...$ and $H_1(N) =
\Z/79\Z$.  The manifold $N$ has a genus-2 Heegaard splitting, and is
the 2-fold branched cover of the 3-bridge knot in Figure~\ref{knot}.
One of the reasons that $N$ was so difficult is that $\pi_1(N)$ has very few
low-index subgroups (the smallest index is 13).  In the end, a search
using \texttt{Magma} \cite{MAGMA}, turned up a subgroup of index 14
which has positive betti number.  It is very hard to enumerate all
finite-index subgroups for an index as large as 14, roughly because
the size of $S_n$ is $n!$; finding this index 14 subgroup took 2 days
of computer time.

While $\pi_1(N)$ has few subgroups of low index, it does have a
reasonable number of simple quotients, and might be a good place to
look for a co-final sequence of covers which fail to have positive
betti number.   The manifold $N$ is non-Haken, but it contains a
essential lamination (and thus a genuine lamination
\cite{Calegari:promoting}).  Arithmetically, it is quite a complicated
manifold---\Snap\ \cite{Snap} computes that the trace field has a
minimal polynomial $p(x)$ whose degree is 51 and largest coefficient
is about $4 \times 10^7$.  The coefficients of $p$ are, starting with
the constant term: 

{
\begin{sloppypar}
\raggedright \small $1$, $24$, $223$, $929$, $909$, $-6163$, $-20232$,
$-2935$, $79745$, $121259$, $-57077$, $-428280$, $-507427$, $689749$,
$2245466$, $-519994$, $-5455251$, $355551$, $9513149$, $-1958013$,
$-12213255$,
 $7478063$, $10535124$, $-17696676$, $-4109720$,
$30159462$, $-2803266$, $-39076707$, $5291640$, $39199917$,
$-3032906$, $-30650313$, $-365203$, $18711624$, $1997701$, $-8892931$,
$-1776338$, $3259601$, $951237$, $-903591$, 
$-352258$, $182336$, $93101$, $-24677$, $-17396$, $1748$, $2197$, $33$, $-169$, $-17$, $6$, $1$.

\end{sloppypar}
}

\subsection{Overlap with known results}

The manifolds we examined have little overlap with those covered by
the known results about the Virtual Haken Conjecture.  The only
general results are those of Cooper and Long \cite{CooperLong99,
  CooperLongSSSSS} building on work of Freedman and Freedman
\cite{FreedmanFreedman}. These are Dehn surgery results---they say
that many ``large'' Dehn fillings on a 1-cusped hyperbolic 3-manifold
are virtually Haken. Because ``large'' Dehn fillings usually have short
geodesics, the Cooper-Long results probably apply to very few, if any,
of the census manifolds.

\subsection{Limitations}

It's possible the behavior we found might not be true in general
because the census manifolds are non-generic in a couple ways.  First,
they all have fundamental groups with presentations with at most 3
generators.  About 75\% have 2-generator presentations.  For these
manifolds, it seems that (at least most of the time) the number of
generators and the Heegaard genus coincide.  So most of these
manifolds have Heegaard genus 2 or 3.

\begin{figure}[ht!]
\begin{center}
\includegraphics[scale=0.6]{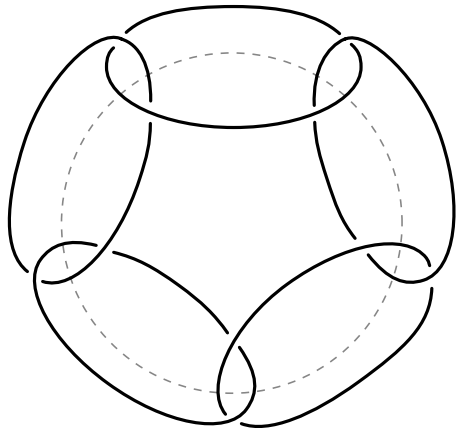}
\end{center}
\caption{The minimally-twisted 5-chain link.}\label{5chain}
\end{figure}

Moreover Callahan, Hodgson, and Weeks (unpublished) showed that almost
all of the census manifolds are Dehn surgeries on a single 5-component
link, the minimally twisted 5-chain shown in Figure~\ref{5chain}.
Let $L$ be this link and $M = S^3 \setminus N(L)$ be its exterior.  The link
$L$ is invariant under rotation of $\pi$ about the dotted grey axis.
The induced involution of $M$ acts on each torus in $\bdry M$ by the
elliptic involution.  Thus the involution of $M$ extends to an
orientation preserving involution of every Dehn filling of $M$.  So
almost all of the census manifolds have an orientation preserving
involution where the fixed point set is a link and underlying space of
the quotient is $S^3$.  While any manifold which has a genus-2
Heegaard splitting has such an involution \cite{BirmanHilden}, this
says that the other 25\% of the census manifolds are also special.  The
presence of such an involution has proven useful in the past.  For
instance, it implies that the manifold is geometrizable.  So it's
possible that our computations only reflect the situation for
manifolds of this type.

\begin{figure}[ht!]
\begin{center}
\includegraphics[scale=0.7]{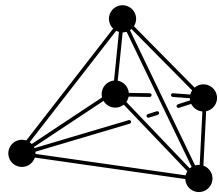}
\end{center}
\caption{The pentacle.}\label{pentacle}
\end{figure}

The 5-chain $L$ is a truly beautiful link, and it's worth describing
some of its properties here.    The orbifold $N$ which is $M$ modulo this
involution is easy to describe.  Take the triangulation $\T$ of $S^3$
gotten by thinking of $S^3$ as the boundary of the 4-simplex.  The
1-skeleton of $\T$ is called the \emph{pentacle}, see
Figure~\ref{pentacle}. 
 If we take $S^3$ minus an
open ball about each vertex in $\T$, and label what's left of each
edge of the pentacle by $\Z/2\Z$, we get exactly the orbifold $N$!  

We can put a hyperbolic structure on $N$ and thus $M$ by making each
tetrahedron in $\T$ a regular ideal tetrahedron.  Thus the volume of
$M$ is $10 v_3 = 10.149416064...$, and further $M$ is arithmetic and
commensurable with the Bianchi group $\PSL{2}{\mathcal{O}(
  \sqrt{-3})}$.  The symmetric group $S_5$ acts on the 4-simplex by
permuting the vertices, inducing an action of $S_5$ on $N$.  This
action is exactly the group of isometries of $N$.  The isometry group
of $M$ is $S_5 \times \Z/2\Z$, where the $\Z/2\Z$ is the rotation about the
axis.

The manifold $M$ fibers over the circle, and in fact every face of the
Thurston norm ball is fibered.  Here's an explicit way to see that $N$
fibers over the interval $I$ with mirrored endpoints (this fibration
lifts to a fibration of $M$ over $S^1$).  Take any Hamiltonian cycle
in the 1-skeleton of $\T$.  The complementary edges also form a
Hamiltonian cycle.  Split the fat vertices of $\T$ (the cusps of $N$)
in the obvious way in space so that these two cycles become the
unlink, with cusps stretched between them.  Then the special fibers
over the $\Z/2\Z$ endpoints of $I$ are two pentagons, spanning the two
Hamiltonian cycles.  The other fibers are 5-punctured spheres.

\section{Techniques for computing homology}\label{techniques}

Given a finite index subgroup $H$ of a finitely presented group $G$, a
simplified version of the Reidemeister-Schreier method produces a matrix
$A$ with integers entries whose cokernel is the abelianization of $H$.
Computing this matrix is not very time-consuming.  The hard part of
computing the rank of the abelianization of $H$ is finding the rank of
$A$. Computing the rank of a matrix is $O(n^3)$ if field operations are
constant time.  We need to compute the rank over $\Q$ so the time needed
is somewhat more than that (see Section~\ref{comp_rank}).   The side
lengths of $A$ are usually about $n = [G:H]$, which at $O(n^3)$ is
prohibitive for many of the covers that we looked at (the largest
covering group we needed was $\PSL{2}{\F_{101}}$, whose order is
515,100).

So one wants to keep the degree of the cover, or really the size of the
matrix involved, as small as possible.  One way to do this, first used
in this context by Holt and Plesken~\cite{HoltPlesken92}, is the
following application of the representation theory of finite groups.
Suppose $H$ is a finite index subgroup of $G$.  Assume that $H$ is
normal, so the corresponding cover is regular.   Set $Q = G/H$ and let 
$f\maps G \to Q$ be the quotient map.  The group $Q$ acts on the
homology of the cover $H_1(H, \C)$, giving a representation of $Q$ on
the vector space $H_1(H, \C)$.  Another description of $H_1(H, \C)$ is
that it is the homology with twisted coefficients $H_1(G, \C Q)$.  As a
$Q$-module, $\C Q$ decomposes as $\C Q = V_1^{n_1} \oplus V_2^{n_2}
\oplus \dots \oplus V_k^{n_k}$ where the $V_i$ are simple $Q$-modules
and $\dim V_i = n_i$. 
So 
\[
H_1(H) = H_1(G, \C Q) = H_1(G, V_1)^{n_1} \oplus H_1(G, V_2)^{n_2}
\oplus \dots \oplus H_1(G, V_k)^{n_k}.
\]
Since the dimensions of the $V_i$ are usually much less than the order
of $Q$, the matrices involved in computing $H_1(G, V_i)$ are much
smaller than the one you would get by applying Reidemeister-Schreier
to the subgroup $H$.  For instance, $\PSL{2}{\F_p}$ has order about
$(1/2) p^3$, but every $V_i$ has dimension about $p$.  If we want to
show that $H_1(H, \C)$ is non-zero, we just have to compute that a
single $H_1(G, V_i)$ is non-zero.

There are a couple of difficulties in computing $H_1(G, V_i)$.  First,
to do the computation rigorously, we need to compute not over $\C$ but
over a finite extension of $\Q$.  Now there is a field $k$ so that $k
Q$ splits over $k$ the same way as $\C Q$ splits over $\C$.  However,
the matrices we need to compute $H_1(G, V_i)$ will have entries in
$k$, whereas the matrix given to us by Reidemeister-Schreier has
integer entries.  If $A$ is a matrix with entries in $k$, to compute
its rank over $\Q$ one can form an associated $\Q$-matrix $B$ by
embedding $k$ as a subalgebra of $\GL{n}{\Q}$ where $n$ is $[k : \Q]$
(see e.g.\  \cite{PleskenSouvignier98}).  The rank of $B$ can then be
computed using one the techniques for integer matrices.  However, the
size of $B$ is the size of $A$ times $[k : \Q]$, so this eats up part
of the apparent advantage to computing just the $H_1(G, V_i)$.

The other problem is that we may not know what the irreducible
representations of $Q$ are, especially if we don't know much about
$Q$.  While computing the character table of a finite group is a
well-studied problem, the problem of finding the actual
representations is harder and not one of the things that \GAP\ or
other standard programs can do.  Even when the representations of $Q$
are explicitly known (e.g.~$Q = \PSL{2}{\F_p}$), it can be
time-consuming to tell the computer how to construct the
representations.  For more on computing the actual representations see
\cite{Dixon70, PleskenSouvignier97}.

We used the following modified approach which avoids
the two difficulties just mentioned, while still reducing the size of
the matrices considerably.  Suppose we are given normal subgroup $H$
and we want to determine if $H_1(H, \C)$ is non-zero.  Suppose $U$ is
a subgroup of $Q$.  Note $U$ is not assumed to be normal.  The
permutation representation of $Q$ on $\C[Q/U]$ desums into
 irreducible representations, say $\C[Q/U] = V_1^{e_1} \oplus V_2^{e_2}
\oplus \dots \oplus V_k^{e_k}$.  Let $K = f^{-1}(U)$, a finite index
subgroup of $G$ containing $H$.  Then
\[
H_1(K) =  H_1(G,\C[Q/U]) =  H_1(G, V_1)^{e_1} \oplus H_1(G, V_2)^{e_2}  \oplus \dots \oplus
H_1(G, V_k)^{e_k}.
\]
Suppose that $U$ is chosen so that every irreducible representation
appears in $\C[Q/U]$, that is, every $e_i > 0$.  Then we see that
$H_1(H)$ is non-zero if and only if $H_1(K)$ is.  As long as $U$ is
non-trivial, the index $[G:K] = [Q:U]$ is smaller than $[G:H] = \#Q$,
so computing $H_1(K)$ is easier that computing $H_1(H)$.  Returning
to the example of $\PSL{2}{\F_p}$, there is such a $U$ of index about
$p^2$, whereas the order of $\PSL{2}{\F_p}$ is about $p^3/2$.  Looking
at a matrix with side $O(p^2)$ is a big savings over one of side
$O(p^3)$.

Moreover, finding such a $U$ given $Q$ is easy.  First compute the
character table of $Q$ and the conjugacy classes of subgroups of $Q$
(these are both well-studied problems).  For each subgroup $U$ of $Q$
compute the character $\chi_U$ of the permutation representation of
$Q$ on $\C[Q/U]$.  Expressing $\chi_U$ as a linear combination of the     
irreducible characters tells us exactly what the $e_i$ are.  Running
through the $U$, we can find the subgroup of lowest index where all of
the $e_i > 0$.

When we were searching for positive betti number covers, we used this
method of replacing $H$ with $K = f^{-1}(U)$ and computed the ranks of
the resulting matrices over a finite field $\F_p$.  Once we had found
an $H$ with positive $\F_p$-betti number, we did the following to
check rigorously that $H$ has infinite abelianization.  First, we went
through \emph{all} the subgroups $U$ of $Q$, till we found the $U$ of
smallest index such that $f^{-1}(U)$ has positive $\F_p$-betti number.
For this $U$, we computed the $\Q$-betti number of $f^{-1}(U)$ using
one of the methods described in Section~\ref{comp_rank}.  Doing this
kept the matrices that we needed to compute the $\Q$-rank of small,
and was the key to checking that the covers really had positive
$\Q$-betti number.  For instance, for the $\PSL{2}{\F_{101}}$-cover of
degree 515,100 there was a $U$ so that the intermediate cover
$f^{-1}(U)$ with positive betti number had degree ``only'' 5,050.

It's worth mentioning that the rank over $\Q$ was very rarely different
than that over a small finite field.  Initially, for each manifold we
found a cover where the $\F_{31991}$-betti number was positive. All but
3 of those 10,986 covers had positive $\Q$-betti number.

\section{Computing the rank over $\Q$}\label{comp_rank}

Here, we describe how we computed the $\Q$-rank of the
matrices produced in the last section.  Normally, one thinks of linear
algebra as ``easy'', but standard row-reduction is polynomial time
only if field operations are constant time.  To compute the rank of an
integer matrix $A$ rigorously one has to work over $\Q$.  Here, doing
row reduction causes the size of the fractions involved to
explode.  There are a number of ways to try to avoid this.

The first is to use a clever pivoting strategy to minimize the size
of the fractions involved \cite{HavasMajewski97, HavasMajewski94,
  HavasHoltRees}.  This is the method built into \GAP , and was
what we used for the covers of degree less than 500, which sufficed
for $99.2\%$ of the manifolds.

For all but about 7 of the remaining 94 manifolds, we used a
simplified version of the $p$-adic algorithm of Dixon given in
\cite{Dixon83}.  Over a large finite field $\F_p$, we computed a basis
of the kernel of the matrix.  Then we used ``rational reconstruction'',
a partial inverse to the map $\Q \to \F_p$ to try to lift each of the
$\F_p$-vectors to $\Q$-vectors (see \cite[pg.~139]{Dixon83}).  If we
succeeded, we then checked that the lifted vectors were actually in the
kernel over $\Q$.  

For 7 of the largest covers (degree 1,000--5,000), this simplification
of Dixon's algorithm fails, and we used the program \MAGMA\ 
\cite{MAGMA}, which has a very sophisticated $p$-adic algorithm, to
check the ranks of the matrices involved.

\section{Simple covers}\label{simple}

To gain more insight into this problem, we looked at a range of simple
covers for a randomly selected 1,000 of the census manifolds which have
2-generator fundamental groups.  For these 1,000 manifolds we found all
the covers where the covering group was a non-abelian finite simple
group of order less than 33,000.  For each cover we computed the
homology.  We will describe some interesting patterns we found.

\begin{table}[ht!]
\begin{center}
\begin{tabular}{rrcccc}
\textbf{Quotient}  & \textbf{Order} &  \textbf{Hit} &  \textbf{HavCov} & \textbf{SucRat1} & \textbf{SucRat2} \\
$A_5$     &    60 & 14.0 &  26.9  &  52.0  &  52.9  \\
$L_2(7)$  &   168 & 17.8 &  28.2  &  63.1  &  66.3 \\
$A_6$     &   360 & 21.6 &  31.4  &  68.8  &  68.7\\
$L_2(8)$  &   504 & 15.4 &  21.7  &  71.0  &  72.6\\
$L_2(11)$ &   660 & 24.1 &  32.8  &  73.5  &  71.8\\
$L_2(13)$ &  1092 & 29.4 &  41.1  &  71.5  &  77.8\\
$L_2(17)$ &  2448 & 29.4 &  43.1  &  68.2  &  69.6\\
$A_7$     &  2520 & 41.1 &  45.8  &  89.7  &  90.9\\
$L_2(19)$ &  3420 & 28.2 &  44.4  &  63.5  &  65.7\\
$L_2(16)$ &  4080 & 11.3 &  18.3  &  61.7  &  65.3\\
$L_3(3)$  &  5616 & 19.2 &  28.0  &  68.6  &  76.5\\
$U_3(3)$  &  6048 & 16.4 &  18.0  &  91.1  &  92.8\\
$L_2(23)$ &  6072 & 32.7 &  47.6  &  68.7  &  70.1\\
$L_2(25)$ &  7800 & 24.7 &  33.0  &  74.8  &  75.5\\
$M_{11}$    &  7920 & 14.6 &  17.1  &  85.4  &  88.8\\
$L_2(27)$ &  9828 & 14.2 &  26.6  &  53.4  &  57.1\\
$L_2(29)$ & 12180 & 42.0 &  57.1  &  73.6  &  74.1\\
$L_2(31)$ & 14880 & 38.1 &  56.5  &  67.4  &  70.9\\
$A_8$     & 20160 & 18.7 &  20.7  &  90.3  &  92.3\\
$L_3(4)$  & 20160 & 42.8 &  50.2  &  85.3  &  89.1\\
$L_2(37)$ & 25308 & 24.9 &  54.2  &  45.9  &  50.5\\
$U_4(2)$  & 25920 & 26.6 &  27.8  &  95.7  &  97.5\\
$Sz(8)$   & 29120 & 26.9 &  43.9  &  61.3  &  73.1\\
$L_2(32)$ & 32736 & 12.4 &  17.9  &  69.3  &  72.1\\
 & & & & & \\
\end{tabular}
\end{center}

\caption{
\textbf{Hit} is the percentage of manifolds having a cover 
with this group which has positive betti number.  \textbf{HavCov} is the
 percentage of manifolds having a cover with this group.
\textbf{SucRate1} is the percentage of manifolds having a cover with this
 group which have such a cover with positive betti number. 
 \textbf{SucRate2} is the percentage of covers with this group having positive
 betti number.}
\label{basic_table}
\end{table}

First, look at Table~\ref{basic_table}.  There, the simple groups are
listed by their ATLAS \cite{ATLAS} name (so, for instance, $L_n(q) =
\PSL{n}{\F_q}$), together with basic information about how many covers
there are, and how many have positive betti number.  There is quite a
bit of variation among the different groups.  For instance, only
$11.3\%$ of the manifold groups have $L_2(16)$ quotients but $42.8\%$
have $L_3(4)$ quotients.  Moreover, there are big differences in how
successful the different kinds of covers are at producing homology.
Only half of the $L_2(37)$ covers have positive betti number, but
almost all ($97.5\%$) of the $U_4(2)$ covers do.  There are no obvious
reasons for these patterns (for instance, the success rates don't
correlate strongly with the order of the group). It would be very
interesting to have heuristics which explain them, and we will explore
these issues in \cite{DunfieldThurstonHeuristics}.

\begin{figure}[ht!]
\begin{center}
\includegraphics[scale=1.0]{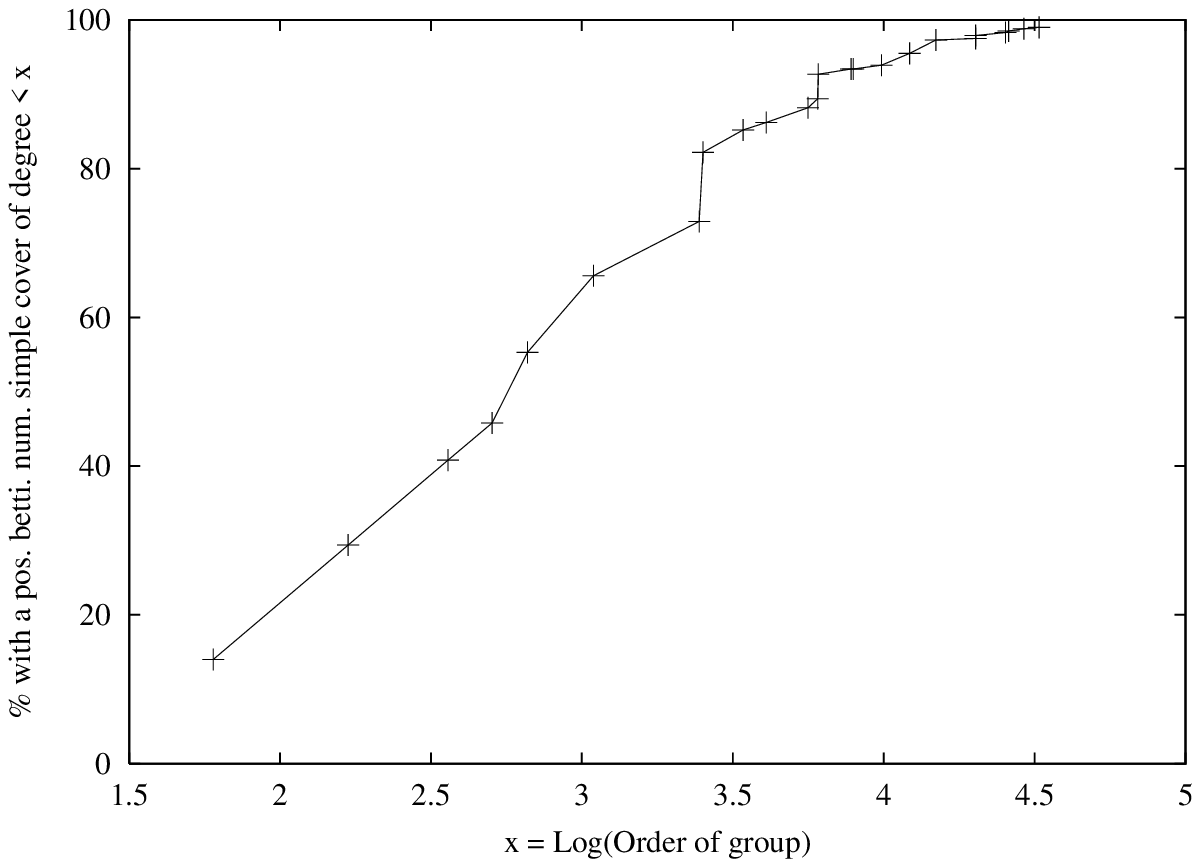}
\end{center}
\caption{
  This graph shows how quickly simple group covers generate homology.
  Each $+$ plotted is the pair $(\log( \#Q(n)) , V(n))$, where the
  $\log$ is base 10.  Thus the leftmost $+$ corresponds to the fact
  that $14\%$ of the manifolds have an $A_5$ cover with positive betti
  number.  The second leftmost $+$ corresponds to the fact that $29\%$
  of the manifolds have either an $A_5$ or an $L_2(7)$ cover with
  positive betti number, and so on.}\label{simple_groups}
\end{figure}

In terms of showing manifolds are virtually Haken, even the least
useful group has a \textbf{Hit} rate greater than $10\%$.  That is,
for any given group at least $10\%$ of the manifolds have a positive
betti number cover with that group.  So unless things are strongly
correlated between different groups, one would expect that every
manifold would have a positive betti number simple cover, and that one
would generally find such a cover quickly.  Let $Q(n)$ denote the
$n^{\mathrm{th}}$ simple group as listed in Table~\ref{basic_table}.
Set $V(n)$ to be the proportion of the manifolds which have a positive
betti number $Q(k)$-cover where $k \leq n$.  We expect that the
increasing function $V(n)$ should rapidly approach $1$ as $n$
increases.  This is born out in Figure~\ref{simple_groups}.

Figure~\ref{simple_groups} shows that the groups behave pretty
independently of each other, although not completely as we will see.  Let
$H(n)$ denote the hit rate for $Q(n)$, that is the proportion of the
manifolds with a $Q(n)$ cover with positive betti number.  If
everything were independent, then one would expect
\[
V(n) \approx V(n-1) + (1-V(n-1)) H(n).
\]
If we let $E(n)$ be the right-hand side above, and compare $E(n)$ with
$V(n)$ we find that $E(n) - V(n)$ is almost always positive.  To judge
the size of this deviation, we look at
\[
\frac{E(n) - V(n)}{1 - V(n-1)} \mtext{which lies in $[-0.007, 0.13]$,}
\]
and which averages $0.022$.  In other words, $V(n) - V(n-1)$ is
usually about $2\%$ smaller as a proportion of the possible increase
than $E(n) - V(n-1)$.

\begin{figure}[ht!]
\begin{center}
\includegraphics[scale=0.9]{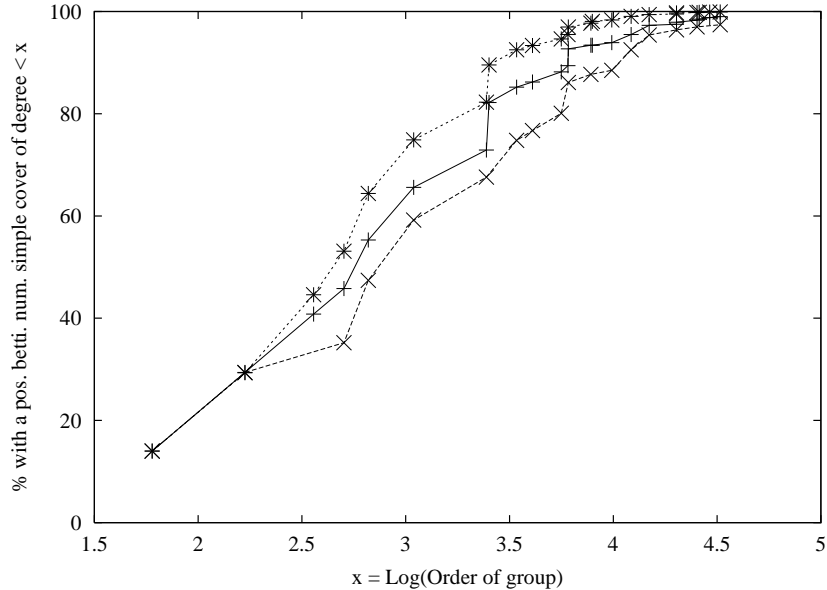}
\end{center}
\caption{The top line plots $(\log(\#Q(n)), V'(n))$, the middle line 
  $(\log(\#Q(n)), V(n))$ (as in Figure~\ref{simple_groups}), and the
  lowest line plots only the groups of the form $L_2(q)$.}\label{simple_compare}
\end{figure}

For a graphical comparison, define $V'(n)$ by the recursion
\[
V'(n) = V'(n-1) + (1-V'(n-1)) H(n),
\]
and compare with $V(n)$ in Figure~\ref{simple_compare}.  

Asymptotically, every non-abelian finite simple group is of the form
$L_2(q)$, and so it's interesting to look at a modified $V(n)$ where we
look only at the $Q(n)$ of this form.  This is also shown in
Figure~\ref{simple_compare}.

\subsection{Amount of homology}  

 \begin{figure}[ht!]
\begin{center}
\psfrag{b}{$\log(R(n))$}
\includegraphics[scale=0.8]{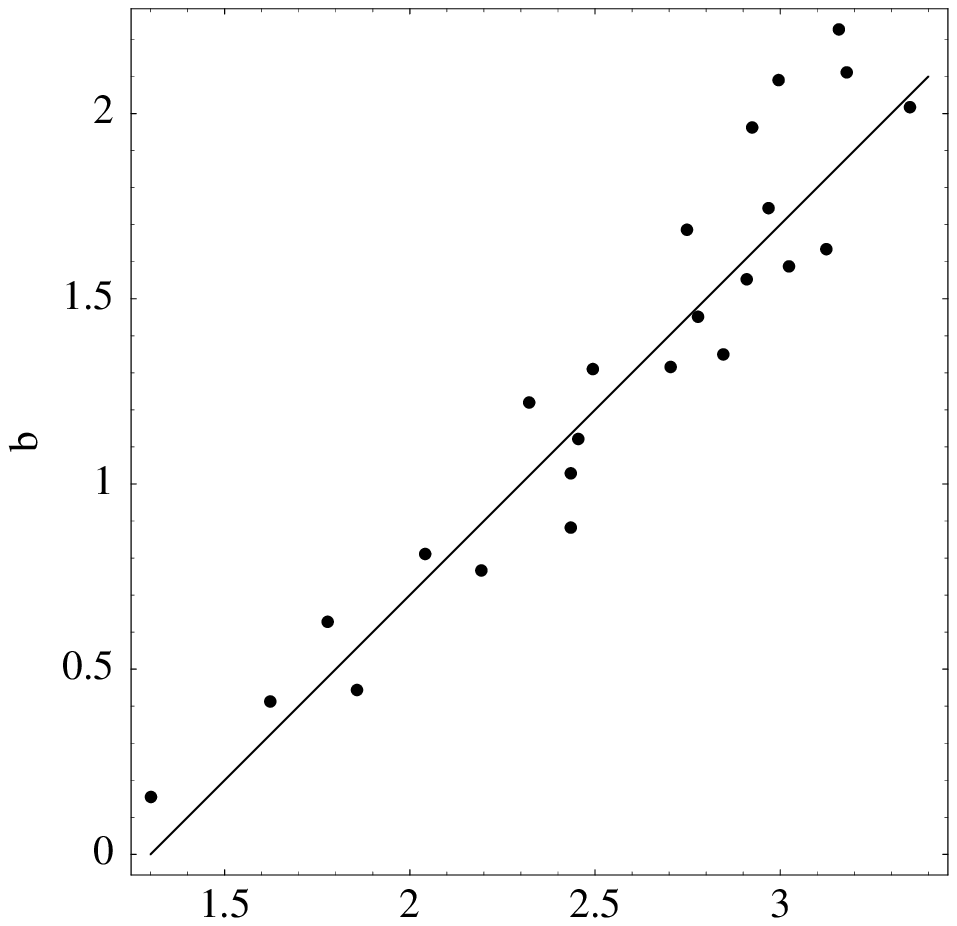}
\end{center}
\vspace{-0.65cm}\centerline{$\log(\#Q(n))$}
\caption{This plot shows the relationship between the expected rank
 and the degree of the cover.  The line shown is $y = x - 1.3$.
  }\label{homology}
\end{figure}

Suppose we look at a simple cover of degree $d$, what is the expected
rank of the homology of the cover?  The data suggests that the
expected rank is linearly proportional to $d$.  For the simple group
$Q(n)$, set $R(n)$ to be the mean of $\beta_1(N)$, where $N$ runs over
all the $Q(n)$ covers of our manifolds (including those where
$\beta_1(N) = 0$).  Figure~\ref{homology} gives a plot of $\log R(n)$
versus $\log(\#Q(n))$.  Also shown is the line $y = x - 1.3$ (which is almost
the least squares fit line $y = 1.018 x - 1.303$).  The data points
follow that line, suggesting that:
\begin{equation}
 \log R(n) \approx \log(\#Q(n)) - 1.3 \mtext{and hence} R(n) \approx
 \frac{\#Q(n)}{20}. \label{hom_rel}
\end{equation}
Now each of the 3-manifold groups we are looking at here are quotients
of the free group on two generators $F_2$.  Let $G$ be fundamental
group of one of our 3-manifolds, say $G = F_2/N$.  Given a
homomorphism $G \to Q(n)$, we can look at the composite homomorphism
$F_2 \to Q(n)$.  Let $H$ be the kernel of $G \to Q(n)$ and $K$ the
kernel of $F_2 \to Q(n)$.  Then the rank of $H_1(K)$ is $\#Q(n) + 1$.
As $H_1(H)$ is a quotient of $H_1(K)$, Equation~\ref{hom_rel} is says
that on average, $5\%$ of $H_1(K)$ survives to $H_1(H)$.

This amount of homology is not a priori forced by the high hit rate
for the $Q(n)$.  For instance, $L_2(p)$ has order $(p^3 - p)/2$ but
has a rational representation of dimension $p$.  Thus it would be
possible for $L_2(p)$ covers to have
\[
\log(R(n)) \approx (1/3) \log(\#G(n))+ C,
\]
even if a large percentage of these covers had positive betti
number.  This data suggests that on a statistical level these
3-manifold groups are trying to behave like the fundamental group of a
2-dimensional orbifold of Euler characteristic $-1/20$.

\sh{Caveats}  The data in Figure~\ref{homology} is not based on the 
full $Q(n)$ covers but on subcovers coming from a fixed subgroup $U(n)
< Q(n)$, chosen as described in Section~\ref{techniques}.  The degree
plotted is the degree of the cover that was used, that is $[Q(n) :
U(n)]$ not the order of $Q(n)$ itself, so the above analysis is still
valid.  Also, throughout Section~\ref{simple} having positive betti
number really means having positive betti number over $\F_{31991}$.
Also, we originally used a list of the Hodgson-Weeks census which had
a few duplicates and so there are actually 12 manifold which appear
twice in our list of 1000 random manifolds.  

\subsection{Homology of particular representations}

As discussed in Section~\ref{techniques}, if we look at a cover with
covering group $Q$, the homology of the cover decomposes into
\[
H_1(G, V_1)^{n_1} \oplus H_1(G, V_2)^{n_2}  \oplus \dots \oplus
H_1(G, V_k)^{n_k},
\]
where $G$ is the fundamental group of the base manifold and the $V_i$
are the irreducible $Q$-modules.  For $Q$ an alternating group, we
looked at this decomposition and found that the ranks of the $H_1(G,
V_i)$ were very strongly positively correlated.  This contrasts with
the relative independence of the ranks of covers with different
$Q(n)$.

\begin{figure}[ht!]
\begin{center}
\includegraphics[scale=0.7]{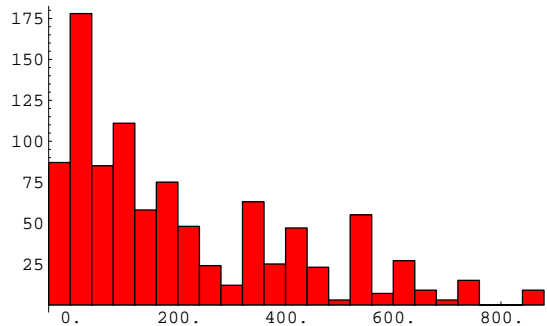}
\end{center}
\caption{Plot showing the distribution of the ranks of the homology of the
  964 covers with group $A_7$.  The $x$-axis is the amount of homology
  and the $y$-axis the number of covers with homology in that range.}
  \label{alt_homology}
\end{figure}

We will describe what happens for $A_7$, the other alternating groups
being similar.  The rational representations of $A_7$ are easy to
describe: they are the restrictions of the irreducible representations
of $S_7$.  They correspond to certain partitions of $7$.
Table~\ref{alt_basic_table} lists the representations and their basic
properties.  Table~\ref{alt_correlations} shows the correlations
between the ranks of the $H_1(G, V_i)$.  Many of the correlations are
larger than $0.5$ and all are bigger than $0$ ($+1$ is perfect
correlation, $-1$ perfect anti-correlation and $0$ the expected
correlation for independent random variables).
Figure~\ref{alt_homology} shows the distribution of the homology of
the covers.

\begin{table}
\begin{center}
\begin{tabular} {rrrr}
Partition  & Dim.~of rep & Success rate & Mean homology \\  
7 & 1 &  2\%  &  0.0 \\
$1,6$ &  6 & 22\% & 1.5 \\
$2,5$ & 14 & 63\% &  19.8 \\
$1,1,5$ & 15  & 64\% & 21.8 \\
$3,4$ & 14 &  41\% &  11.0 \\
$1,2,4$ & 35 & 70\% &  101.6 \\
$1,1,1,4$ & 20  & 61\% & 20.7 \\
$1,3,3$ & 21 & 61\% &  33.9
\end{tabular}
\end{center}

\caption{ The $\Q$-irreducible representations of $A_7$.  Success Rate
 is the percentage of covers where that representation appeared.  
 Mean Homology is the average amount of homology that that representation 
 contributed (the mean homology of an $A_7$ cover was 210.3). }
\label{alt_basic_table}
\end{table}

\begin{table}
\begin{center}
\begin{tabular} {r|cccccccc}
 &  7 & 16 & 25 & 115 & 34 & 124 & 1114 & 133 \\
  \hline
7 &\textbf{1.00} & 0.01 & 0.11 & 0.08 & 0.15 & 0.17 & 0.02 & 0.13 \\
16  & 0.01 & \textbf{1.00} & 0.22 & 0.09 & 0.23 & 0.19 & 0.18 & 0.19 \\
25  & 0.11 & 0.22 & \textbf{1.00} & 0.63 & 0.65 & 0.79 & 0.37 & 0.61  \\
115 & 0.08 & 0.09 & 0.63 & \textbf{1.00} & 0.52 & 0.80 & 0.75 & 0.78 \\
34  & 0.15 & 0.23 & 0.65 & 0.52 & \textbf{1.00} & 0.73 & 0.50 & 0.65 \\
124 & 0.17 & 0.19 & 0.79 & 0.80 & 0.73 & \textbf{1.00} & 0.65 & 0.89 \\
1114& 0.02 & 0.18 & 0.37 & 0.75 & 0.50 & 0.65 & \textbf{1.00} & 0.66 \\
133 & 0.13 & 0.19 & 0.61 & 0.78 & 0.65 & 0.89 & 0.66 & \textbf{1.00} 
\end{tabular}
\end{center}
\vspace{0.2in}
\caption{Table showing the correlations between the ranks of $H_1(G, V_i)$ where the $V_i$ are indexed by the partition of the corresponding representation. } \label{alt_correlations}
\end{table}

\subsection{Correlations between groups}\label{correlation}

In the beginning of Section~\ref{simple} we saw that the two events
\[
\big(\mbox{having a $Q(n)$-cover with $\beta_1 > 0$}, \mbox{having a $Q(m)$-cover with $\beta_1>0$}\big)
\]
were more or less independent of each other, though overall there was
a slight positive correlation which dampened the growth of $V(n)$.  In
the appendix, there is a table giving these
correlations, was well one giving those between the events:
\[
\big(\mbox{having a $Q(n)$-cover},  \mbox{having a $Q(m)$-cover}\big).
\]
Some of these correlations are much larger than one would expect by
chance alone---for instance the correlation between
\[
\big(\mbox{ having a $L_2(7)$-cover with $\beta_1 > 0$}, 
\mbox{ having a $L_2(8)$-cover with $\beta_1 > 0$}\big)
\]
is $0.38$.  Moreover, there are very few negative correlations and
those that exist are quite small.  Overall, the average correlation is
positive as we would expect from Section~\ref{simple}.

One way of trying to understand these correlations is to observe that
almost all of these manifolds are Dehn surgeries on the minimally
twisted $5$-chain.  Let us focus on the simpler question of
correlations between having a cover with group $Q(n)$
and having a cover with group $Q(m)$.  Let $M$ be the complement of
the $5$-chain.  Consider all the homomorphisms $f_k \maps \pi_1 M \to
Q(n)$.  Supposes $X$ is a Dehn filling on $M$ along the five slopes
$(\gamma_1, \gamma_2, \gamma_3, \gamma_4, \gamma_5)$ where $\gamma_i$
is in $\pi_1(\bdry_i M)$. The manifold $X$ has a cover with group
$Q(n)$ if and only if there is an $f_k$ where each $\gamma_i$ lies in
the kernel of $f_k$ restricted to $\pi_1( \bdry_i M)$.  Thus having a
cover with group $Q(n)$ is determined by certain subgroups of the
groups $\pi_1(\bdry_i M) = \Z^2$.  If we consider a different group
$Q(m)$ we get a different family of subgroups of the $\pi_1(\bdry_i
M)$.  If there is a lot of overlap between these two sets of
subgroups, there will be a positive correlation between having a cover
with group $Q(n)$ and having a cover with group $Q(m)$.  If there is
little overlap then there will be a negative correlation.  However,
even looked at this way there seems to be no reason that the
average correlation should be positive.

If we look at the same question for manifolds which are Dehn surgeries
on the figure-8 knot (a simplified version of this setup) there are
many negative correlations and the overall average correlation is 0.
If we look at the question for small surgeries on the Whitehead link,
the overall average correlation is positive and of similar magnitude
of that for the 5-chain.  If we also look at larger surgeries on the
Whitehead link the average correlation drops somewhat.  By changing
the link we get a different pattern of correlations, and so it is
unwise to attach much significance to these numbers.

\section{Further questions}\label{questions}

Here are some interesting further questions related to our
experiment.
\begin{enumerate}
  
\item What happens for 3-manifolds bigger than the ones we looked at?
  Do the patterns we found persist?  It is computationally difficult
  to deal with groups with large numbers of generators, which would
  limit the maximum size of the manifolds considered.  But another
  difficulty is how to find a ``representative'' collection of such
  manifolds.  (Some notions of a ``random 3-manifold'', which help
  with this latter question, will be discussed in
  \cite{DunfieldThurstonHeuristics}).
  
\item How else could the virtually Haken covers we found be used to
  give insight into these conjectures?  For instance, one could try to
  look at the virtual fibration conjecture.  While there is no good
  algorithm for showing that a closed manifold is fibered, one could
  look at the following algebraic stand-in for this question.  If a
  3-manifold fibers over the circle, then one of the coefficients of
  the Alexander polynomial which is on a vertex of the Newton polytope
  is $\pm 1$ (see e.g.~\cite{Dunfield:norms}).  One could compute the
  Alexander polynomial of the covers with virtual positive betti
  number and see how often this occurred.  As many of our covers are
  quite small, computing the Alexander polynomial should be feasible
  in many cases.
  
\item One could use our methods to look at the Virtual Positive Betti
  Number conjecture for lattices in the other rank-1 groups that don't
  have Property T.  This would be particularly interesting for the
  examples of complex hyperbolic manifolds where every congruence
  cover has $\beta_1 = 0$.  These complex hyperbolic manifolds were
  discovered by Rogawski \cite[Thm.~15.3.1]{Rogawski90} and are
  arithmetic.
  
\end{enumerate}

\section{Transferring virtual Haken via Dehn filling}\label{transfer}

In the rest of this paper, we consider the following setup.  Let $M$
be a compact 3-manifold with boundary a torus. The process of
\textit{Dehn filling} creates closed 3-manifolds from $M$ by taking a
solid torus $D^2 \times S^1$ and gluing its boundary to the boundary
of $M$.  The resulting manifolds are parameterized by the isotopy class
of essential simple closed curve in $\bdry M$ which bounds a disc in
the attached solid torus.  If $\alpha$ denotes such a class, called a
\emph{slope}, the corresponding Dehn filling is denoted by
$M(\alpha)$.  Though no orientation of $\alpha$ is needed for Dehn
filling, we will often think of the possible $\alpha$ as being the
primitive elements in $H_1(\bdry M, \Z)$ and so $H_1(\bdry M, \Z)$
parameterizes the possible Dehn fillings.

If you have a general conjecture which you can't prove for all 3-manifolds, a
standard thing to do is to try to prove it for most Dehn fillings on
an arbitrary 3-manifold with torus boundary.  For instance, in
the case of the Geometrization Conjecture there is the following theorem:

\begin{hyperbolic_dehn_surgery}{\rm\cite{ThurstonLectureNotes}}\qua
  Let $M$ be a compact 3-man\-i\-fold with $\bdry M$ a torus.  Suppose
  the  interior of $M$ has a complete hyperbolic metric of finite
  volume.  Then all but finitely many Dehn fillings of $M$ are
  hyperbolic manifolds.
\end{hyperbolic_dehn_surgery}

For the Virtual Haken Conjecture there is the following result of
Cooper and Long.  A properly embedded compact surface $S$ in $M$ is
\emph{essential} if it is incompressible, boundary incompressible, and
not boundary parallel.  Suppose $S$ is an essential surface in $M$.
While $S$ may have several boundary components, they are all parallel
and so have the same slope, called the boundary slope of $S$. If
$\alpha$ and $\beta$ are two slopes, we denote their minimal
intersection number, or \emph{distance}, by $\Delta(\alpha, \beta)$.

\begin{theorem}{\rm (Cooper-Long \cite{CooperLong99})}\qua
  Let $M$ be a compact orientable 3-manifold with torus boundary which
  is hyperbolic.  Suppose $S$ is a non-separating orientable essential
  surface in $M$ with non-empty boundary.  Suppose that $S$ is not the
  fiber in a fibration over $S^1$.  Let $\lambda$ be the boundary
  slope of $S$.  Then there is a constant $N$ such that for all slopes
  $\alpha$ with $\Delta(\alpha, \lambda) \geq N$, the manifold
  $M(\alpha)$ is virtually Haken.
  
  Explicitly, $N = 12 g - 8 + 4 b$ where $g$ is the
  genus of $S$ and $b$ is the number of boundary components.
\end{theorem}
 
This result differs from the Hyperbolic Dehn Surgery Theorem in that
it excludes those fillings lying in an infinite strip in $H_1(\bdry
M)$, instead of only excluding those in a compact set.  Here, we will
prove a Dehn surgery theorem about the Virtual Positive Betti Number
Conjecture, assuming that $M$ has a very simple Dehn filling which
strongly has virtual positive betti number.  Our theorem is a
generalization of the work of Boyer and Zhang \cite{BoyerZhang00},
which we discuss below.
 
The basic idea is this.  Suppose $M$ has a Dehn filling $M(\alpha)$
which has virtual betti number in a very strong way.  By this we mean
that there is a surjection $\pi_1(M(\alpha)) \to \Gamma$ where
$\Gamma$ is a group all of  whose finite index subgroups have lots of
homology.  In our application, $\Gamma$ will be the fundamental group
of a hyperbolic 2-orbifold.  Given some other Dehn filling $M(\beta)$,
we would like to transfer virtual positive betti number from
$M(\alpha)$ to $M(\beta)$.  Look at $\pi_1(M)/\pair{\alpha, \beta}$
which we will call $\pi_1(M(\alpha, \beta))$.  This group is a common
quotient of $\pi_1(M(\alpha))$ and $\pi_1(M(\beta))$.  Choose $\gamma
\in \pi_1(\bdry M)$ so that $\{\alpha, \gamma\}$ is a basis of
$\pi_1(\bdry M)$.  Then $\beta = \alpha^m \gamma^n$.  If we think of
$\pi_1(M(\alpha, \beta))$ as a quotient of $\pi_1(M(\alpha))$ we
have:
\[
\pi_1(M(\alpha, \beta)) = \pi_1(M(\alpha)) /\pair{\beta} = \pi_1(M(\alpha))/\pair{\gamma^n}.
\]
Thus $\pi_1(M(\alpha, \beta))$ surjects onto $\Gamma/\pair{\gamma^n}$,
where here we are confusing $\gamma$ and its image in $\Gamma$.  So
$\pi_1(M)$ surjects onto $\Gamma/\pair{\gamma^n}$.  If $\Gamma$ has rapid
homology growth, one can hope that $\Gamma_n = \Gamma/\pair{\gamma^n}$
still has virtual positive betti number when $n$ is large
enough.  This is plausible because adding a relator which is
a large power often doesn't change the group too much.  If
there is an $N$ so that $\Gamma_n$ has virtual positive betti number
for all $n \geq N$, then $M(\beta)$ has virtual positive betti number
for all $\beta$ with $n = \Delta(\gamma, \alpha) \geq N$.

Our main theorem applies when $M(\alpha)$ is a Seifert fibered space
whose base orbifold is hyperbolic:
\begin{theorem}\label{SF_filling}
  Let $M$ be a compact 3-manifold with boundary a torus.  Suppose
  $M(\alpha)$ is Seifert fibered with base orbifold $\Sigma$
  hyperbolic.  Assume also that the image of $\pi_1(\bdry M)$ under
  the induced map $\pi_1(M) \to \pi_1(\Sigma)$ contains no non-trivial
  element of finite order.  Then there exists an $N$ so that
  $M(\beta)$ has virtual positive betti number whenever
  $\Delta(\alpha, \beta) \geq N$.
  
  If $\Sigma$ is not a sphere with 3 cone points, then $N$ can be taken
  to be $7$.
\end{theorem}
In light of the above discussion, if we consider the homomorphism
$\pi_1(M(\alpha)) \to \pi_1(\Sigma) = \Gamma$,
Theorem~\ref{SF_filling} follows immediately from:
\begin{theorem}\label{orbifold_quotient}
  Let $\Sigma$ be a closed hyperbolic 2-orbifold without mirrors, and
  $\Gamma$ be its fundamental group.  Let $\gamma \in \Gamma$ be a
  element of infinite order.  Then there exists an $N$ such that for
  all $n \geq N$ the group
  \[
  \Gamma_n = \Gamma/\pair{\gamma^n}
  \]
  has virtual positive betti number. In fact, $\Gamma_n$ has a finite
  index subgroup which surjects onto a free group of rank 2.
  
  If $\Sigma$ is not a 2-sphere with 3 cone points, then $N =
  \max\{1/\abs{1 + \chi(\Sigma)}, 3\}$.  In this case, $N$ is at most
  $7$.
\end{theorem}

In applying Theorem~\ref{SF_filling}, the technical condition that the
image of $\pi_1(\bdry M)$ not contain an element of finite order holds in
many cases.  For instance, Theorem~\ref{SF_filling} implies the
following theorem about Dehn surgeries on the Whitehead link.  Let $W$
the exterior of the Whitehead link.  Given a slope $\alpha$ on the
first boundary component of $W$, we denote by $W(\alpha)$ the manifold
with one torus boundary component obtained by filling along $\alpha$.
\begin{whiteforward}
  Let $W$ be  the exterior of the Whitehead link. Then for all but
  finitely many slopes $\alpha$, the manifold $M = W(\alpha)$ has the
  following property:  All but finitely many Dehn fillings of $M$ have
  virtual positive betti number.  
\end{whiteforward}
In fact, our proof of this theorem excludes only 28 possible slopes
$\alpha$ (see Section~\ref{Whitehead}).  The complements of the twist
knots in $S^3$ are exactly the $W(1/n)$ for $n \in \Z$.
Theorem~\ref{whiteheadthm} applies to all of the slopes $1/n$ except for
$n \in \{0, 1\}$ which correspond to the unknot and the trefoil.  Thus we
have:
\begin{corollary}
  Let $K$ be a twist knot in $S^3$ which is not the unknot or the
  trefoil.  Then all but finitely many Dehn surgeries on $K$ have
  virtual positive betti number.
\end{corollary}
For the simplest hyperbolic knot, the figure-8, we can use a
quantitative version of Theorem~\ref{orbifold_quotient} due to Holt and
Plesken \cite{HoltPlesken92} which applies in this special case.
We will show:
\begin{theorem}
  Every non-trivial Dehn surgery on the figure-8 knot in $S^3$ has
  virtual positive betti number.
\end{theorem}

As we mentioned, Theorem~\ref{SF_filling} generalizes the work of
Boyer and Zhang \cite{BoyerZhang00}.  They restricted to the case
where the base orbifold was not a 2-sphere with 3 cone points. In
particular, they proved:
\begin{theorem}{\rm\cite{BoyerZhang00}}\label{BoyerZhang}\qua
  Let $M$ have boundary a torus.  Suppose $M(\alpha)$ is Seifert
  fibered with a hyperbolic base orbifold $\Sigma$ which is not a
  2-sphere with 3 cone points. Assume also that $M$ is small, that
  is, contains no closed essential surface.  Then $M(\beta)$ has
  virtual positive betti number whenever $\Delta( \alpha, \beta) \geq
  7$.
\end{theorem}
The condition that $M$ is small is a natural one as if $M$ contains an
closed essential surface, then there is a $\alpha$ so that $M(\beta)$
is actually Haken if $\Delta(\alpha, \beta) > 1$ \cite{CGLS, Wu91}.

Boyer and Zhang's point of view is different than ours, in that they
do not set out a restricted version of
Theorem~\ref{orbifold_quotient}.  While the basic approach of both
proofs comes from \cite{BaumslagMorganShalen}, Boyer and Zhang's proof
of Theorem~\ref{BoyerZhang} also uses the Culler-Shalen theory of
$\SL{2}{\C}$-character varieties and surfaces arising from ideal
points.  From our point of view this is not needed, and
Theorem~\ref{BoyerZhang} follows easily from Theorem~\ref{SF_filling}
(see the end of Section~\ref{one_rel_quo} for a proof).

In Section~\ref{dehn_questions}, we discuss possible generalizations
of Theorem~\ref{SF_filling} to other types of fillings.  In a very
special case, we use toroidal Dehn fillings to show
(Theorem~\ref{sister_thm}) that every Dehn filling of the sister of
the figure-8 complement satisfies the Virtual Positive Betti Number
Conjecture.

\section{One-relator quotients of 2-orbifold groups}\label{one_rel_quo}

This section is devoted to the proof of
Theorem~\ref{orbifold_quotient}.  The basic ideas go back to
\cite{BaumslagMorganShalen} which proves the analogous result for
$\Gamma = \Z/p * \Z/q$.  Fine, Roehl, and Rosenberger proved
Theorem~\ref{orbifold_quotient} in many, but not all, cases where
$\Sigma$ is not a 2-sphere with 3 cone points
\cite{FineRoehlRosenberger, FineRosenberger99}.  In the case $\Sigma =
S^2(a_1, a_2, a_3)$, Darren Long and Alan Reid suggested the proof
given below, and Matt Baker provided invaluable help with the number
theoretic details.

\begin{proof}[Proof of Theorem \ref{orbifold_quotient}]
  
Let $\Sigma_n$ be the 2-complex with marked cone points consisting of
$\Sigma$ together with a disc $D$ with a cone point of order $n$,
where the boundary of $D$ is attached to $\Sigma$ along a curve
representing $\gamma$.  Thus $\Gamma_n = \pi_1(\Sigma_n)$.  Now the
Euler characteristic of $\Sigma_n$ is $\chi(\Sigma) + 1/n$, which is
negative if $n > 1/\abs{\chi(\Sigma)}$.  From now on, assume that $n >
1/\abs{\chi(\Sigma)}$.  Suppose $\Gamma_n$ contains a subgroup
$\Gamma_n'$ of finite index such that if $\alpha$ is a small loop
about a cone point then $\alpha \not\in \Gamma_n'$.  For instance,
this is the case if $\Gamma_n'$ is torsion free.  Let $\Sigma_n'$ be
the corresponding cover of $\Sigma_n$, so $\Gamma_n' =
\pi_1(\Sigma_n')$.  Then $\Sigma_n'$ is a 2-complex without any cone
points.  Since $\Sigma_n'$ has negative Euler characteristic and there
is no homology in dimensions greater than two, we must have
$H_1(\Sigma_n', \Q) \neq 0$.  Thus $\Gamma_n$ has virtual positive betti
number.
  
One can show more: Let $d$ be the degree of the cover $\Sigma_n' \to
\Sigma_n$.  The complex $\Sigma_n'$ is a smooth hyperbolic surface $S$
with $d/n$ discs attached.  From this description it is easy to check
that $\Gamma_n'$ has a presentation where
\begin{equation*}
\begin{split}
(\mbox{\# of generators}) - (\mbox{\# of relations}) &= (\abs{\chi(S)}+ 1) - \frac{d}{n}  \\
&= 1 + d\left( \abs{\chi(\Sigma)} - \frac{1}{n}
\right) \geq 2.
\end{split}
\end{equation*}
By a theorem of Baumslag and Pride \cite{BaumslagPride78}, the group
$\Gamma_n'$ has a finite-index subgroup which surjects onto $\Z * \Z$.
  
So it remains to produce the subgroups $\Gamma_n'$. First, we discuss
the case where $\Sigma$ is not a sphere with 3 cone points.  A
homomorphism $f \maps \Gamma \to Q$ is said to preserve torsion if for
every torsion element $\alpha$ in $\Gamma$ the order of $f(\alpha)$ is
equal to the order of $\alpha$.  (Recall that the torsion elements of
$\Gamma$ are exactly the loops around cone points.)  The key is to
show:
\begin{lemma}\label{non_triangle_case}
  Suppose $\Sigma$ is not a 2-sphere with 3 cone points, and that
  $\gamma \in \Gamma$ has infinite order.   Given any $n > 2$, there
  exists a homomorphism $\rho \maps \Gamma \to \PSL{2}{\C}$ such that
  $\rho$ preserves torsion and $\rho(\gamma)$ has order $n$.
\end{lemma}
Suppose we have $\rho$ as in the lemma, which we will regard as a
homomorphism from $\Gamma_n$ to $\PSL{2}{\C}$.  By Selberg's lemma,
the group $\rho(\Gamma)$ has a finite index subgroup $\Lambda$ which
is torsion free.  We can then take $\Gamma_n'$ to be
$\rho^{-1}(\Lambda)$.  Because the lemma only requires that $n>2$ and
the preceding argument required that $n>1 /\abs{\chi(\Sigma)}$, in
this case we can take the $N$ in the statement of
Theorem~\ref{orbifold_quotient} to be $\max\{3, 1 +
1/\abs{\chi(\Sigma)}\}$.  A case check, done in \cite{BoyerZhang00},
shows that $N$ is at most $7$.  As we will see, the proof of
Lemma~\ref{non_triangle_case} is relatively easy and involves
deforming Fuchsian representations $\Gamma \to \Isom(\H^2)$ to find
$\rho$.
   
The harder case is when $\Sigma$ is a 2-sphere with 3 cone points,
which we denote $S^2(a_1, a_2, a_3)$.  Here the fundamental group
$\Gamma$ can be presented as
\[
  \spandef{x_1, x_2, x_3 }{ x_1^{a_1} = x_2^{a_2} = x_3^{a_3} = x_1
    x_2 x_3 = 1}.
\]
Geometrically, $x_i$ is a loop around the $i$\textsuperscript{th} cone
point.  We will show:
\begin{lemma}\label{triangle_case}
  Let $\Gamma = \pi_1(S^2(a_1, a_2, a_3))$ where $1/a_1 + 1/a_2 +
  1/a_3 < 1$.  Given an element $\gamma \in \Gamma$ of infinite order,
  there exists an $N$ such that for all $n \geq N$ the group $\Gamma$
  has a finite quotient where the images of $(x_1, x_2, x_3, \gamma)$
  have orders exactly $(a_1, a_2, a_3, n)$  respectively.
\end{lemma}
With this Lemma, we can take $\Gamma'_n$ to be the kernel of the given
finite quotient.  The proof of Lemma~\ref{triangle_case} involves
using congruence quotients of $\Gamma$ and a some number theory.
Unfortunately, unlike the previous case, the proof of
Lemma~\ref{triangle_case} gives no explicit bound on $N$.
  
In any event, we've established Theorem~\ref{orbifold_quotient} modulo
Lemmas~\ref{non_triangle_case} and \ref{triangle_case}.
\end{proof}
  
The rest of this section is devoted to proving the two lemmas.  

\begin{proof}[Proof of Lemma~\ref{non_triangle_case}]
  
  Because $\Sigma$ is not a 2-sphere with 3 cone points, the
  Teichm{\"u}ller space of $\Sigma$ is positive dimensional.  Thus there
  are many representations of $\Gamma$ into $\Isom(\H^2)$.  We can
  embed $\Isom(\H^2)$ into $\Isom^+(\H^3) = \PSL{2}{\C}$ as the
  stabilizer of a geodesic plane.  We will then deform these Fuchsian
  representations to produce $\rho$.
  
  Pick a simple closed curve $\beta$ which intersects $\gamma$
  essentially.  There are two cases depending on whether a
  neighborhood of $\beta$ is an annulus or a M{\"o}bius band.  
  
  Suppose the neighborhood is an annulus.  First, let's consider the
  case where $\beta$ separates $\Sigma$ into 2 pieces.  In this case
  $\Gamma$ is a free product with amalgamation $A *_{\pair{\beta}} B$.
  Let $\rho_1 \maps \Gamma \to \PSL{2}{\C}$ be one of the Fuchsian
  representations.  Conjugate $\rho_1$ so that $\rho_1(\beta)$ is
  diagonal.  Then $\rho_1(\beta)$ commutes with the matrices
  \[
  C_t = \left(\begin{array}{cc} t & 0 \\
                                0 & t^{-1}
        \end{array}\right)  \quad  \mbox{for $t$ in $\C^\times$}.
  \]  
  For $t$ in $\C^\times$, let $\rho_t$ be the representation of
  $\Gamma$ whose restriction to $A$ is $\rho_1$ and whose restriction
  to $B$ is $C_t \rho_1 C_t^{-1}$.  Consider the function $f \maps
  \C^\times \to \C$ which sends $t$ to $\tr^2(\rho_t(\gamma))$.  It is
  easy to see that $f$ is a rational function of $t$ by expressing
  $\gamma$ as a word in elements of $A$ and $B$.  We claim that $f$ is
  non-constant.  First, suppose that neither of the two components of
  $\Sigma \setminus \beta$ is a disc with two cone points of order 2.
  In this case, $\beta$ can be taken to be a geodesic loop.  If we
  restrict $t$ to $\R$ then the family $\{ \rho_t \}$ corresponds to
  twisting around $\beta$ in the Fenchel-Nielsen coordinates on
  $\Teich(\Sigma)$.  As $\gamma$ intersects $\beta$ essentially, the
  length of $\gamma$ changes under this twisting and so $f$ is
  non-constant.  From this same point of view, we see that that $f$
  has poles at $0$ and $\infty$.  If one of the pieces of $\Sigma
  \setminus \beta$ is a disc with two cone points of order 2, then
  $\beta$ naturally shrinks not to a closed geodesic, but to a
  geodesic arc joining the two cone points.  There is still a
  Fenchel-Nielsen twist about $\beta$, and so  we have the same
  observations about $f$ in this case (think of $\Sigma$ being
  obtained from a surface with a geodesic boundary component by
  pinching the boundary to a interval).

  Since the rational function $f$ has poles at $\{ 0, \infty\}$, we
  have $f(\C^\times) = \C$.  So given $n > 1$, we can choose $t \in
  \C^\times$ so that $\tr^2(\rho_t(\gamma)) = (\zeta_{2n} +
  \zeta_{2n}^{-1})^2$ where $\zeta_{2n} = e^{\pi i /n}$.  Then
  $\rho_t(\gamma)$ has order $n$.  Moreover, $\rho_t$ preserves
  torsion because $\rho_1$ does, and so we have finished the proof of
  the lemma when $\beta$ is separating and has an annulus
  neighborhood.  If $\beta$ has an annulus neighborhood and is
  non-separating, the proof is identical except that $\Gamma$ is an
  HNN-extension instead of a free product with amalgamation.
  
  Now we consider the case where the neighborhood of $\beta$ is a
  M{\"o}bius band.  The difference here is that you can't twist a
  hyperbolic structure of $\Sigma$ along $\beta$.  To see this, think
  of constructing $\Sigma$ from a surface with geodesic boundary where
  the boundary is identified by the antipodal map to form $\beta$.
  Instead, we will deform the length of $\beta$ in $\Teich(\Sigma)$.
  Here we will need the hypothesis that $n > 2$, as you can see by
  looking at $\RP^2(3,5)$ with $\gamma$ a simple closed geodesic which
  has a M{\"o}bius band neighborhood.  The only quotient of
  $\pi_1(\RP^2(3,5))$ where $\gamma$ has order $2$ is $\Z/2$ and this
  doesn't preserve torsion.
  
  The underlying surface of $\Sigma$ is non-orientable.  We can assume
  that $\Sigma$ has at least one cone point since every non-orientable
  surface covers such an orbifold.  Pick an arc $a$ joining $\beta$ to
  a cone point $p$.  Let $A$ be a closed neighborhood of $\beta \cup
  a$.  The set $A$ is a M{\"o}bius band with a cone point.  Let $B$ be
  the closure of $\Sigma \setminus A$.  Let $\alpha$ be the boundary
  of $A$.  A small neighborhood of $\alpha$ is an annulus, so if
  $\gamma$ intersects $\alpha$ essentially, we can replace $\beta$
  with $\alpha$ and use the argument above.  So from now on, we can
  assume that $\gamma$ lies in $A$.  Let $\psi \maps \Gamma \to
  \PSL{2}{\C}$ be a Fuchsian representation.  Suppose we construct a
  representation $\rho \maps \pi_1(A) \to \PSL{2}{\C}$ so that $\rho$
  preserves torsion, $\rho(\gamma)$ has order $n$, and
  $\tr^2(\rho(\alpha)) = \tr^2(\psi(\alpha))$.  Then as $\Gamma =
  \pi_1(A) *_{\pair{\alpha}} \pi_1(B)$ and $\rho$ and $\psi$ are
  conjugate on $\pair{\alpha}$, we can glue $\rho$ and $\psi$
  restricted to $\pi_1(B)$ together to get the required representation
  of $\Gamma$.
  
 Thus we have reduced everything to a question about certain
  representations of $\pi_1(A)$.  The group $\pi_1(A)$ is generated by
  $\alpha$ and $\beta$.  Choosing orientations correctly, a small loop
  about the cone point $p$ is $\delta = \beta^2 \alpha$.   If $p$ has
  order $r$, then $\pi_1(A)$ has the presentation
  \[
  \spandef{\alpha, \beta, \delta}{\delta = \beta^2 \alpha, \delta^r = 1}.
  \]
  
  Given any representation $\phi$ of $\pi_1(A)$, we will fix lifts of
  $\phi(\alpha)$ and $\phi(\beta)$ to $\SL{2}{\C}$.  Having done this,
  any word $w$ in $\alpha$ and $\beta$ has a canonical lift of
  $\phi(w)$ to $\SL{2}{\C}$.  We will abuse notion and denote this
  lift by $\phi(w)$ as well.  In this way, we can treat $\phi$ as
  though it was a representation into $\SL{2}{\C}$ so that, for
  instance, the trace of $\phi(w)$ is defined.
  
  Define a 1-parameter family of representations $\rho_t$ for $t \in
  \C^\times$ as follows.  Set
  \[
  \psi(\beta) = \left(\begin{array}{cc} 
      0 & 1 \\
      -1 & t 
    \end{array}\right),
  \mtext{and}
  \psi(\alpha) = \left(\begin{array}{cc}
      e &   s \\
      0 &  e^{-1}
    \end{array}\right)
  \]
  where $e + e^{-1} = \tr(\psi(\alpha))$ and $s = \frac{1}{t}(e^{-1} t^2 - (e + e^{-1}) - \tr(\psi(\delta))$.  This gives a representation
  of $\pi_1(A)$ because $s$ was chosen so that $\tr(\rho_t(\delta)) =
  \tr(\psi(\delta))$ and so $\rho_t(\delta)$ also has order $r$ in
  $\PSL{2}{\C}$.    
  
  Let $\Teich(A)$ denote hyperbolic structures on $A$ with geodesic
  boundary where the length of the boundary is fixed to be that of the
  Fuchsian representation $\psi$.  This Teichm{\"u}ller space is $\R$
  with the single Fenchel-Nielsen coordinate being the length of
  $\beta$.  Note that any irreducible representation of $\pi_1(A)$ is
  conjugate to some $\rho_t$, and so each point in $\Teich(A)$ yields
  a Fuchsian representation $\rho_t$.  As $\beta$ gets short in
  $\Teich(A)$, the curve $\gamma$ gets long.  Thus if we set $f =
  \tr(\rho_t(\gamma))$, then $f$ is a non-constant Laurent
  polynomial in $t$.
  
  Let $v = \zeta_{2n} + \zeta^{-1}_{2n}$.  To finish the proof of the
  lemma, all we need to do is find a $t \in \C^\times$ so that $f(t)^2
  = v^2$.  As a map from the Riemann sphere to itself, $f$ is onto and
  there are $t_1$ and $t_2$ in $\widehat{\C}$ so that $f(t_1) = v$ and
  $f(t_2) = -v$.  As $n > 2$, $v$ is not $0$ and so $t_1$ and $t_2$
  are distinct.  As $f$ is non-constant and finite on $\C^\times$, it
  has a pole at at least one of $0$ and $\infty$.  Therefore, at least
  one of $t_1$ and $t_2$ is in $\C^\times$ and we are done.
\end{proof}

\begin{proof}[Proof of Lemma~\ref{triangle_case}] 
  
  The group $\Gamma$ is naturally a subgroup of $\PSL{2}{\R}$.  Set
  $b_i = 2 a_i$.   Let $X_i$ be the matrix in $\PSL{2}{\R}$
  corresponding to the generator $x_i$.  As $X_i$ has order $a_i$, it
  follows that $\tr(X_i) = \pm (\zeta_{b_i} + \zeta_{b_i}^{-1})$ where
  $\zeta_{b_i}$ is some primitive $b_i$\textsuperscript{th} root of
  unity.  Any irreducible 2-generator subgroup of $\PSL{2}{\C}$ is
  determined by its traces on the generators and their product, and so
  we can conjugate $\Gamma$ in $\PSL{2}{\C}$ so the $X_i$ are:
  \[ X_1 = \left(\begin{array}{cc} 
      0 & 1 \\
      -1 & \zeta_{b_1} + \zeta_{b_1}^{-1}
    \end{array}\right),
  X_2 = \left(\begin{array}{cc}
      \zeta_{b_2} +\zeta_{b_2}^{-1} & -\zeta_{b_3}  \\
      \zeta_{b_3}^{-1} & 0
    \end{array}\right),
  \ \mbox{and}\  X_3 = (X_1 X_2)^{-1}.
  \]
  Henceforth we will identify $\Gamma$ with its image.  The entries of
  the $X_i$ lie in $\Q(\zeta_{b_1}, \zeta_{b_2}, \zeta_{b_3})$, and
  moreover are integral, so $\Gamma$ is contained in the subgroup
  $\PSL{2}{{\O}(\Q(\zeta_{b_1}, \zeta_{b_2}, \zeta_{b_3}))}$.  Let $G$
  be a matrix in $\PSL{2}{\C}$ representing $\gamma$.  Let $a$ be one
  of the eigenvalues of $G$.  Note that $a$ is an algebraic integer,
  in fact a unit, because it satisfies the equation $a^2 - (\tr G) a +
  1$ and $\tr{G}$ is integral.  Let $K$ be the field $\Q(\zeta_{b_1},
  \zeta_{b_2}, \zeta_{b_3}, a)$.  From now on, we will consider
  $\Gamma$ as a subgroup of $\PSL{2}{{\O}(K)}$.  We will construct the
  required quotients of $\Gamma$ from congruence quotients of
  $\PSL{2}{{\O}(K)}$.  Suppose $\p$ is a prime ideal of ${\O}(K)$. Setting
  $k = {\O}(K)/\p$, we have the finite quotient of $\Gamma$ given
  by
  \[
  \Gamma \to \PSL{2}{{\O}(K)} \to \PSL{2}{k}.
  \]
  What conditions do we need so that $(x_1, x_2, x_3, \gamma)$ have
  the right orders in $\PSL{2}{k}$?  Well, the eigenvalues of $X_i$
  are $\pm \{ \zeta_{b_i}, \zeta_{b_i}^{-1} \}$, so as long as
  $\bar{\zeta}_{b_i}$ has order $b_i$ in $k^\times$, the matrix
  $\bar{X}_i$ in $\PSL{2}{k}$ also has order $b_i$.  Similarly, if we
  set $m = 2 n$, then $\bar{G}$ in $\PSL{2}{k}$ has order $n$ if
  $\bar{a}$ has order $m$ in $k^\times$.  Thus the following claim
  will complete the proof of the lemma:
  \begin{claim}
    There exists an $N$ such that for all $n \geq N$ there is a prime
    ideal $\p$ such that if $k = {\O}(K)/\p$ then the images of
    $(\zeta_{b_1}, \zeta_{b_2}, \zeta_{b_3}, a)$ in $k^\times$ have
    orders $(b_1, b_2, b_3, m)$.
  \end{claim}
  
  Let's prove the claim.  The idea is to show that $a^m - 1$ is not a
  unit in ${\O}(K)$ for large $m$, and then just take $\p$ to be a prime
  ideal dividing $a^m - 1$.  We have to be careful, though, that
  $(\bar{\zeta}_{b_1}, \bar{\zeta}_{b_2}, \bar{\zeta}_{b_3}, \bar{a})$
  don't end up with lower orders that expected in $k^\times$.
      
  A prime ideal is called \emph{primitive} if it divides $a^m - 1$ and
  does not divide $a^r - 1$ for all $r < m$.  Postnikova and Schinzel
  proved the following theorem:
  \begin{theorem}{\rm\cite{Schinzel74, PostnikovaSchinzel}}\label{num_theory}\qua
    Suppose that $a$ is an algebraic integer which is not a root of
    unity.  There there is an $N$ such that for all $n \geq N$ the
    integer $a^n - 1$ has a primitive divisor.
\end{theorem}
The proof of Theorem~\ref{num_theory} relies on deep theorems of
Gelfond and A.~Baker on the approximation by rationals of logarithms
of algebraic numbers.
      
Because $\gamma$ has infinite order, we know that $a$ is not a root of
unity.  Thus Theorem~\ref{num_theory} applies, and let $N$ be as in
the statement.  By increasing $N$ if necessary, we can ensure that the
primitive divisor $\p$ given Theorem~\ref{num_theory} does not
divide any element of the finite set
\[
  R = \setdef{ \zeta_{b_i}^r  - 1}{ 1 \leq r < b_i}.
\]
Thus for all $m \geq N$, we have a prime ideal $\p$ which divides $a^m
- 1$ but does not divide $a^r - 1$ for $r < m$.  Thus $\bar{a}$ has
order $m$ in $k^\times$.  As $\p$ does not divide any element of $R$,
the element $\bar{\zeta}_{b_i}$ has order $b_i$ in $k^\times$.  This
proves the claim and thus the lemma.
\end{proof}

It would be nice to have given a proof of Lemma~\ref{triangle_case}
which gave an explicit bound on $N$.  The number theory used gives
``an effectively computable constant'' for $N$, but doesn't actually
compute it.  Perhaps there are other proofs of
Lemma~\ref{triangle_case} more like that of
Lemma~\ref{non_triangle_case}.  While $\pi_1(S^2(a_1, a_2, a_3))$ has
only a finite number of representations into $\PSL{2}{\C}$, if one
looks at representations into larger groups there are deformation
spaces where you could hope to play the same game.  For instance, if
one embeds $\H^2$ as a totally geodesic subspace in complex hyperbolic
space $\CH^2$, then a Fuchsian representation deforms to a one real
parameter family in $\Isom^+(\CH^2) = \mathrm{PU}(2,1)$.  One could
instead consider deformations in the space of real-projective
structures, which gives rise to homomorphisms to $\PGL{3}{\R}$
\cite{ChoiGoldman2001}.  In general, the structure of the space
representations of $\pi_1(S^2(a_1, a_2, a_3)) \to \SL{n}{\C}$ is closely
related to the Deligne-Simpson problem \cite{Simpson90}.

We end this section by deducing Boyer and Zhang's original
Theorem~\ref{BoyerZhang} from Theorem~\ref{SF_filling}.

\begin{proof}[Proof of Theorem~\ref{BoyerZhang}]
  Let $M$ be a manifold with torus boundary which is small.  Suppose
  that $M(\alpha)$ is Seifert fibered with hyperbolic base orbifold
  $\Sigma$ which is not sphere with 3 cone points.  We need to check
  that Theorem~\ref{SF_filling} applies.  Let $\beta$ be a curve so
  that $\{ \alpha, \beta \}$ is a basis for $\pi_1(\bdry M)$.  It
  suffices to show the image of $\beta$ does not have finite order in
  $\Gamma = \pi_1(\Sigma)$.  Suppose not.  Then there are infinitely
  many Dehn fillings $M(\gamma_i)$ of $M$ where $\pi_1(M(\gamma_i))$
  surjects onto $\Gamma$.  The orbifold $\Sigma$ contains an essential
  simple closed curve which isn't a loop around a cone point.
  Therefore, $\Gamma$ has non-trivial splitting as a graph of groups
  and so acts non-trivially on a simplicial tree.  Then each
  $\pi_1(M(\gamma_i))$ act non-trivially on a tree and so
  $M(\gamma_i)$ contains an essential surface.  As infinitely many
  fillings contain essential surfaces, a theorem of Hatcher
  \cite{Hatcher82} implies that $M$ contains a closed essential
  surface.  This is contradicts that $M$ is small.  So the image of
  $\beta$ has infinite order and we are done.
\end{proof}

\section{Surgeries on the Whitehead link}\label{Whitehead}

\begin{figure}[ht!]
\begin{center}
  \includegraphics[scale=0.65]{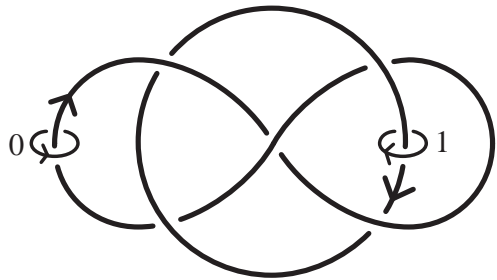}
  \caption{The Whitehead link, showing our orientation conventions
    for the meridians and longitudes.}\label{whitehead_link}
\end{center}
\end{figure}

Consider the Whitehead link pictured in Figure~\ref{whitehead_link}.
Let $W$ be its exterior.
We will denote the two boundary components of $W$ by $\bdry_0 W$ and
$\bdry_1 W$.  For each $\bdry_i W$, we fix a meridian-longitude basis
$\{ \mu_i, \lambda_i \}$ with the orientations shown in the figure.  With
respect to one of these bases, we will write boundary slopes as
rational numbers, where $p \mu + q \lambda$ corresponds to $p/q$.  We will
denote Dehn filling of both boundary components of $W$ by $W(p_0/q_0;
p_1/q_1)$.  Dehn filling on a single component of $W$ will be denoted
$W(p_0/q_0; \unfilled)$ and $W(\unfilled ; p_1/q_1)$.  As $W(p/q;
\unfilled)$ is homeomorphic to $W(\unfilled ; p/q)$, we will sometimes
denote this manifold by $W(p/q)$.  With our conventions, $W(1)$ is the
trefoil complement, and $W(-1)$ is the figure-8 complement.  The
manifold $W(p/q)$ is hyperbolic except when $p/q$ is in $\{ \infty, 0, 1, 2,
3, 4 \}$.  The point of this section is to show:
\begin{theorem}\label{whiteheadthm}
  Let $W$ be the complement of the Whitehead link.  For any slope
  $p/q$ which is not in $E = \{\infty$, $0$, $1$, $2$, $3$, $4$, $5$,
  $5/2$, $6$, $7/1$, $7/2$, $8$, $8/3$, $9/2$, $10/3$, $11/2$, $11/3$,
  $13/3$, $13/4$, $14/3$, $15/4$, $16/3$, $16/5$, $17/5$, $18/5$,
  $19/4$, $24/5$, $24/7\}$ the manifold $W(\alpha)$ has the property
  that all but finitely many Dehn fillings have virtual positive betti
  number.
\end{theorem}

\begin{proof}
  The proof goes by showing that except for $p/q$ in $E$, the manifold
  $W(p/q)$ has at least 2 distinct Dehn fillings which are Seifert
  fibered and to which Theorem~\ref{SF_filling} applies.  The reason
  that $W(p/q)$ has so many Seifert fibered fillings is because the
  manifolds $W(1)$, $W(2)$, and $W(3)$ are all Seifert fibered with
  base orbifold a disc with two cone points.  In particular, the base
  orbifolds are $D^2(2,3)$, $D^2(2,4)$, and $D^2(3,3)$ respectively.
  Therefore, all but one Dehn surgery $W(1;p/q)$ on $W(1)$ is Seifert
  fibered with base orbifold a sphere with 3 cone points.  Similarly
  for $W(2)$ and $W(3)$.  In fact, you can check that
  \begin{itemize}
  \item $W(1; p/q)$ Seifert fibers over $S^2(2,3, \abs{p - 6 q})$ if
    $p/q \neq 6$.
  \item $W(2; p/q)$ Seifert fibers over $S^2(2,4, \abs{p - 4 q})$ if
    $p/q \neq 4$.
  \item $W(3; p/q)$ Seifert fibers over $S^2(3,3, \abs{p - 3 q})$ if
    $p/q \neq 3$.
  \end{itemize}
  
  Now fix a slope $p/q$, and consider the manifold $M = W(\unfilled
  ;p/q)$.  We want to know when we can apply Theorem~\ref{SF_filling}
  to $M(1)$, $M(2)$, or $M(3)$.  First, we need the base orbifold to
  be hyperbolic, i.e.~that the reciprocals of the orders of the cone points
  sum to less than 1.  This leads to the conditions:
  \begin{equation}\begin{split}\label{pq_cond}
  & \mbox{For $M(1)$ that $\abs{p - 6 q} > 6$.} \\ 
  & \mbox{For $M(2)$ that $\abs{p - 4 q} > 4$.} \\
  & \mbox{For $M(3)$ that $\abs{p - 3 q} > 3$.}
  \end{split}
  \end{equation}
  We claim that as long as the base orbifold is hyperbolic then
  Theorem~\ref{SF_filling} applies.  Consider the map $\pi_1(M) \to \Gamma$
  where $\Gamma$ is the fundamental group of one of the base orbifolds.
  Let $\mu$ in $\bdry M$ be the meridian coming from our meridian $\mu_0$
  of $W$.  Since $\mu$ intersects any of the slopes $1, 2, 3$ once, its
  image in $\Gamma$ generates the image of $\pi_1(\bdry M)$.  Thus we just
  need to check that the image of $\mu$ is an element of infinite order
  in $\Gamma$.  One can work out what the image in $\Gamma$ is explicitly
  (most easily by with the help of \SnapPea\ \cite{SnapPea}):
  \begin{equation}\begin{split}\label{quo_pres}
 \mbox{For $M(1)$,} &  \mbox{ $\mu \mapsto aba^{-1}b^{-1}$ where}  \\
    &\Gamma = \spandef{a,b}{ a^2 = b^3 = (ab)^{p - 6 q} = 1}. \\
  \mbox{ For $M(2)$,}& \mbox{ $\mu \mapsto ab^2$ where} \\ 
   & \Gamma = \spandef{a,b}{ a^2 = b^4 = (ab)^{p - 4 q} = 1}. \\
 \mbox{ For $M(3)$,}&  \mbox{ $\mu \mapsto a b^{-1}$ where} \\
    &\Gamma = \spandef{a,b}{ a^3 = b^3 = (ab)^{p - 3 q} = 1}. 
  \end{split}\end{equation}
  It remains to check that the images of $\mu$ above always have
  infinite order in $\Gamma$.  This is intuitively clear for looking
  at loops which represent these elements.  The suspicious reader can
  check that this is really the case by using, say, the solution to
  the word problem for Coxeter groups \cite[{\S} II.3]{Brown89}.
  
  Thus, Theorem~\ref{SF_filling} applies whenever one of the
  conditions in (\ref{pq_cond}) holds.  If $p/q$ is such that two of
  (\ref{pq_cond}) hold, then all but finitely many Dehn surgeries on
  $M$ have virtual positive betti number.  The set in $H_1(\bdry M,
  \R) = \R^2$ where any one of the conditions fails is an infinite
  strip.  So the set where a fixed pair of them fail is compact,
  namely a parallelogram.  Hence, outside a union of 3 parallelograms,
  at least two of the conditions hold.  These 3 parallelograms are all
  contained in the square where $\abs{p}, \abs{q} \leq 100$.  To
  complete the proof of the theorem, one checks all the slopes in that
  square to find those where fewer that two of (\ref{pq_cond})
  hold.
\end{proof}

For most of the slopes in $E$, one of (\ref{pq_cond}) holds, and so
one still has a partial result.  The slopes where none of the
conditions in (\ref{pq_cond}) hold are 
\[
\{ \infty, 0, 1, 2, 3, 4, 5, 6, 7/2, 9/2\}.
\]  
One interesting manifold among these exceptions is
the sister of the figure-8 complement $W(5)$.  We will consider that
manifold in detail in Section~\ref{sister}.

\section{The figure-8 knot}

Here we prove:

\begin{theorem} \label{fig8thm}
  Every non-trivial Dehn surgery on the figure-8 knot has virtual
  positive betti number.
\end{theorem}

\begin{proof}
  Let $M$ be the figure-8 complement.  As the figure-8 knot is
  isotopic to its mirror image, the Dehn filling $M(p/q)$ is
  homeomorphic to $M(-p/q)$.  Now, if $W$ is the Whitehead complement
  as in the last section, $M = W(-1)$.  Hence $M$ has at least 6
  interesting Seifert fibered surgeries namely $M(\pm 1)$, $M(\pm 2)$
  and $M(\pm 3)$.  In (\ref{quo_pres}), we saw exactly which orbifold
  quotients $\Gamma/\pair{\mu^n}$ arise when we try our method of
  transferring virtual positive betti number.  By a minor miracle,
  Holt and Plesken have looked at exactly these quotients and shown:
  \begin{theorem}{\rm\cite{HoltPlesken92}}\qua Let 
    \begin{align*}
    \Gamma^1_n &= \spandef{a,b}{ a^2 = b^3 = (ab)^7 =
      (aba^{-1}b^{-1})^n = 1} ,\\
    \Gamma^2_n &= \spandef{a,b}{ a^2 = b^4 = (ab)^5 = (ab^2)^n} , \mbox{and}\\ 
    \Gamma^3_n &= \spandef{a,b}{ a^3 = b^3 = (ab)^4 = (ab^{-1})^n = 1}.
    \end{align*}
    These groups have virtual positive betti number if $n \geq 11$ for
    $\Gamma^1_n$ and $n \geq 6$ for $\Gamma^2_n$ and $\Gamma^3_n$.
  \end{theorem}
  Thus $M(\alpha)$ has virtual positive betti number if any of the
  following hold: 
  \[
  \Delta(\alpha, \pm 1) \geq 11, \Delta(\alpha, \pm2) \geq 6, \mtext{or} \Delta(\alpha, \pm 3) \geq 6.
  \]
  It's easy to check that the only slopes $\alpha$ for which none of these
  hold are $\{\infty, 0, \pm 1, \pm 2 \}$.  Since $H_1(M(0)) = \Z$ and the
  Seifert fibered manifolds $M(\pm 1)$ and $M(\pm 2)$ have virtual
  positive betti number, we've proved the theorem.
\end{proof}

\section{Other groups of the form $\Gamma/\pair{\gamma^n}$ and further questions}
\label{dehn_questions}

As we have seen, groups of the form $\Gamma/\pair{\gamma^n}$, where
$\Gamma$ is a Fuchsian group, are very useful for studying the Virtual
Haken Conjecture via Dehn filling.  So it is natural to ask: what
other types of $\Gamma$ give similar results?  In this
section, we consider $\Gamma$ which are free products with
amalgamation of finite groups.  The key source here is Lubotzky's
paper \cite{Lubotzky96}, which gives a number of applications of these
groups to the Virtual Positive Betti Number Conjecture.

For convenience, we will only discuss free products with amalgamation,
but there are analogous statements for HNN extensions.  Let $\Gamma =
A *_C B$ be an amalgam of finite groups where $C$ is a proper subgroup
of $A$ and $B$.  The group $\Gamma$ acts on a tree $T$ with finite
point stabilizers.  By \cite[{\S} II.2.6]{SerreTrees}, $\Gamma$ has a
finite index subgroup $\Lambda$ which acts freely on $T$.  The
subgroup $\Lambda$ has to be free, and so $\Gamma$ is virtually free.
It is not hard to show that if one of $[A : C]$ and $[B : C]$ is $\geq
3$ then $\Gamma$ is virtually a free group of rank $\geq 2$
\cite[Lemma 2.2]{Lubotzky96}.  From now on, we will assume $[A : C]
\geq 3$.  Because $\Gamma$ is virtually free, it is natural to hope
that the answer to the following question is yes:
\begin{question}  
  Let $\Gamma$ be an amalgam of finite groups, and fix $\gamma \in
  \Gamma$ of infinite order.  Does there exist an $N$ such that for
  all $n \geq N$, the group $\Gamma_n = \Gamma/\pair{\gamma^n}$ has
  virtual positive betti number?
\end{question}
Note that by Gromov, there is an $N$ such that $\Gamma_n$ is a
non-elementary word hyperbolic group for all $n \geq N$.

Now consider these groups in the context of Dehn filling.  Suppose $M$
is a manifold with torus boundary, and suppose $\alpha$ is a slope
where $\pi_1(M(\alpha))$ surjects onto $\Gamma$, an amalgam of finite
groups.  Choose $\gamma$ in $\pi_1(\bdry M)$ so that $\{ \alpha,
\gamma \}$ form a basis.  The proof of Theorem~\ref{BoyerZhang} shows
that if $M$ does not contain a closed incompressible surface, then the
image of $\gamma$ in $\Gamma$ has infinite order.

There are candidate $\alpha$ where one expects that $\pi_1(M(\alpha))$ will
surject onto an amalgam of finite groups.  Suppose that $N = M(\alpha)$
contains a separating incompressible surface $S$.  Then $\pi_1(N)$
splits as $\pi_1(N_1) *_{\pi_1(S)} \pi_1(N_2)$, where the $N_i$ are the
components of $N \setminus S$.  Recall that $\pi_1(S)$ is said to
\emph{separable} in $\pi_1(N)$ if it is closed in the profinite topology
on $\pi_1(N)$.  Lubotzky showed \cite[Prop.~4.2]{Lubotzky96} that if
$\pi_1(S)$ is separable then there is a homomorphism from $\pi_1(N)$ to an
amalgam of finite groups $\Gamma$, which respects the amalgam structure.
Provided that $S$ is not a semi-fiber (that is, the $N_i$ are not both
$I$-bundles), then $\Gamma = A *_C B$ can be chosen so that $[A:C] \geq 3$.

In general, we will say that $\pi_1(S)$ is \emph{weakly separable}
when there is such an amalgam preserving map from $\pi_1(N)$ to an
amalgam of finite groups.  A priori, this is weaker than $\pi_1(S)$
being closed in $\pi_1(N)$, which is in turn weaker than $\pi_1(N)$
being subgroup separable (aka LERF).

Note that if $\pi_1(S)$ is weakly separable, then $N$ has virtual
positive betti number as $\pi_1(N)$ virtually maps onto a free group.
If $N$ is hyperbolic, it seems quite possible that the fundamental
group of an embedded surface is always weakly separable.  If this is
the case, there is no difference between being virtually Haken and
having virtual positive betti number.  Subgroup separability
properties for 3-manifold groups have been difficult to prove even in
special cases.  Weak separability also seems quite difficult to show
even though the surface $S$ is embedded.

Let $M$ be a manifold with torus boundary which is hyperbolic.  Assume
that $M$ does not contain a closed incompressible surface.  Then there
are always at least two Dehn fillings of $M$ which contain an
incompressible surface \cite{CullerShalen84, CGLS}.  If embedded
surface subgroups are weakly separable, we would expect that for most
$M$, there are at least two slopes where $\pi_1(M(\alpha))$ surjects
onto an amalgam of finite groups.  One has to say ``most'' here
because $M(\alpha)$ might be a (semi-)fiber or the Poincar{\'e}
conjecture might fail.  This makes it plausible that, regardless of
the truth of the virtual Haken conjecture in general, for a fixed $M$
all but finitely many Dehn fillings of $M$ have virtual positive betti
number.  In this context, it is worth mentioning the result of
Cooper-Long \cite{CooperLongSSSSS} which says that for any such
hyperbolic $M$ all but finitely many of the Dehn fillings contain a
surface group.  If fundamental groups of hyperbolic manifolds are
subgroup separable, then this result would also imply that all but
finitely many fillings of $M$ have virtual positive betti number.

One case where weak separability is known is when $N = M(\alpha)$ is
irreducible and the incompressible surface $S$ in $N$ is a torus.  Then
$N$ is Haken and, by geometrization, $\pi_1(N)$ is residually finite.
Using this it's not too hard to show that $\pi_1(S)$ is a separable
subgroup.  So in this case $\pi_1(N)$ maps to a amalgam of finite
groups.  In the next section, we will use these ideas in this special case to
show that all of the Dehn filings on the sister of the figure-8
complement satisfy the Virtual Haken Conjecture.
 
\section{The sister of the figure-8 complement}\label{sister}
 
Let $M$ be the sister of the figure-8 complement.  The manifold $M$ is
the punctured torus bundle where the monodromy has trace $-3$, and is also
the surgery on the Whitehead link $W(5)$.  We will use the basis
$(\mu, \lambda)$ of $\pi_1(\bdry M)$ coming from the standard basis on
$W$.  We will show:
\begin{theorem}\label{sister_thm}
  Let $M$ be the sister of the figure-8 complement.  Then every Dehn
  filling of $M$ which has infinite fundamental group has virtual
  positive betti number.
\end{theorem}

\begin{proof}
  The manifold $M$ has a self-homeomorphism which acts on $\pi_1(\bdry
  M)$ via $(\mu, \lambda) \mapsto (\mu + \lambda, -\lambda)$.  Let $N$
  be the filling $M(4) \iso M(4/3)$.  The manifold $N$ contains a
  separating incompressible torus.  It turns out that this torus
  splits $N$ into a Seifert fibered space with base orbifold
  $D^2(2,3)$ and a twisted interval bundle over the Klein bottle.
  Rather than describe the details of this splitting, we will simply
  exhibit the final homomorphism from $\pi_1(N)$ onto an amalgam of
  finite groups. In fact, $\pi_1(N)$ surjects onto $\Gamma = S_3
  *_{C_2} C_4$ where $C_n$ is a cyclic group of order $n$.

According to \SnapPea, the group $\pi_1(N)$ has
presentation:
\[
\spandef{a, b}{ ab^2ab^{-1}a^3b^{-1} =  ab^2a^{-2}b^2 = 1} 
\]
where $\mu \in \pi(\bdry M)$ becomes $a b$ in $\pi_1(N)$.  If
we add the relators $a^3 = b^4 = 1$ to the presentation of $\pi_1(N)$,
we get a surjection from $\pi_1(N)$ onto
\[
\Gamma  = \spandef{a,b}{a^3 = b^4 = (ab^2)^2 = 1}.
\]
As $S_3$ has presentation $\spandef{x,y}{x^3 = y^2 = (xy)^2=1}$, we see
that $\Gamma$ is $S_3 *_{C_2} C_4$ where the first factor is generated
by $\{a, b^2\}$ and the second by $b$.

We will need:
\begin{lemma}\label{amalgam}
  Let $\Gamma$ be $S_3 *_{C_2} C_4$ and let $\gamma \in \Gamma$ be
  $ab$.  The group
  \[
   \Gamma_n  = \Gamma/\pair{\gamma^n}
  \]
  has virtual positive betti number for all $n \geq 10$.  For $n < 10$, the group
  $\Gamma_n$ is finite.
\end{lemma}

Assuming the lemma, the theorem follows easily.  Given a slope
$\alpha$ in $\pi_1(\bdry M)$, if either $\Delta(\alpha, 4) \geq 10$ or
$\Delta(\alpha, 4/3) \geq 10$ then $M(\alpha)$ has virtual positive
betti number.  The only $\alpha$ which satisfy neither condition are
$E = \{ 0$, $-1$, $\infty$, $1$, $1/2$, $2$, $3$, $3/2$, $4$, $4/3$,
$5/2$, $5/3$, $7/3$, $7/4\}.$ One can check that the fillings along
these slopes either have finite $\pi_1$ or have virtual positive betti
number (the 6 hyperbolic fillings in $E$ are all among the
census manifolds which we showed have virtual positive betti number in
the earlier sections).

Now we will prove the lemma.
\begin{proof}[Proof of Lemma~\ref{amalgam}]

As in  the case of a Fuchsian group the key is to show:

\begin{claim}\label{bigclaim}
  Let $n \geq 12$.  Then there is a homomorphism $f$ from $\Gamma$ to
  a finite group $Q$ where $f$ is injective on the amalgam factors
  $S_3$ and $C_4$ and where $\gamma$ has order $n$.  
\end{claim}

Assuming this claim, we will prove the theorem for $n \geq 12$.  The
Euler characteristic (in the sense of Wall \cite{Wall61}) of $\Gamma$
is $1/6 + 1/4 - 1/2 = -1/12$.  Let $K$ be the kernel of $f$.  The
subgroup $K$ is free, and from its Euler characteristic we see that
it has rank $1 + \#Q/12$.  Let $K'$ be the kernel of the induced
homomorphism from $\Gamma_n \to Q$.  Then $H_1(K', \Z)$ is obtained
from $H_1(K, \Z)$ by adding $\#Q/n$ relators.  As $n \geq 12$, this
implies that $H_1(K', \Z)$ is infinite and $\Gamma_n$ has virtual
positive betti number. 

To prove the rest of the theorem, one can check that
$\Gamma_{10}$ and $\Gamma_{11}$ have homomorphisms into $S_{12}$ and
$\PSL{2}{\F_{23}}$ respectively whose kernels have infinite $H_1$.
Using coset enumeration, it is easy to check that $\Gamma_n$ is
finite for $n < 10$.

Now we establish the claim.  For each $n$, we will inductively build a
permutation representation $f \maps \Gamma \to S_n$ where $f(\gamma)$
has order $n$.  We will say that $f \maps \Gamma \to S_n$ is
\emph{special} if it is faithful on the amalgam factors, $f( \gamma )$
is an $n$-cycle, and $f(b)$ fixes $n$.  If $f$ satisfies these
conditions except for $f(b)$ fixing $n$, we will say that $f$ is
\emph{almost special}.  Our induction tool is:
\begin{claim}\label{subclaim}
  Suppose that $f$ is a special representation of $\Gamma$ into $S_n$.
  Then there exists a special representation of $\Gamma$ into $S_{n +
    6}$.  Also, there exists an almost special representation of $\Gamma$
  into $S_{n +7}$.
\end{claim}
To see this, let $f$ be a special representation.  First, we construct
the representation into $S_{n+6}$.  Let
\[
L = \{1, 2,\dots, n \} \cup \{p_1, p_2, p_3, p_4, p_5, p_6 \}.
\]
We will find a special representation into $S_L$.  Let $g \maps \Gamma
\to S_{\{n, p_1, \dots, p_6\}}$ be the special representation given by
\[
g(a) = ( p_1 p_2 p_3) (p_4 p_5 p_6) \mtext{and} g(b) = (n p_1) (p_2 p_4 p_3 p_5).
\]
It's easy to check (using that $f(a)$ commutes with $g(b^2)$, etc.)
that $h(a) = f(a) g(a)$ and $h(b) = f(b) g(b)$ induces a homomorphism
$h \maps \Gamma \to S_L$.  Moreover, $h( a b ) = f(a) g(a) f(b) g(b) =
f(a)f(b) g(a) g(b) = f(ab) g(ab)$.  Thus $h$ is the product of an
$n$-cycle and a $7$-cycle which overlap only in $n$, and so is a $n +
6$ cycle.  So $h$ is special.

To construct the almost-special representation, do the same thing,
where $g$ replaced is now defined by
\[
g(a) = ( p_1 p_2 p_3) (p_4 p_5 p_6)  \mtext{and} g(b) = (n p_1) (p_2 p_4 p_3 p_5)(p_6 p_7).
\]
This establishes the inductive Claim~\ref{subclaim}.  
 
Using the induction, to prove Claim~\ref{bigclaim} it suffices to show
that there are special representations for $n = 6, 7, 15, 17$,  and
that there is an almost-special representation for $n = 16$.  These are
\begin{align*}
n = 6  \quad &  a \mapsto ( 1,2,3)( 4,5,6)\\
              & b \mapsto ( 2, 4, 3, 5) \\
n = 7  \quad & a \mapsto ( 2, 3, 4)( 5, 6, 7) \\
              & b \mapsto  ( 1, 2)( 3, 5, 4, 6) \\
n = 15 \quad & a \mapsto ( 2, 3, 4)( 5, 7, 9)( 6, 8,11)(12,13,15) \\
                         & b \mapsto  ( 1, 2)( 3, 5, 4, 6)( 7,10,11,14)( 8,12, 9,13)\\
n = 16 \quad &  a \mapsto ( 2, 3, 4)( 5, 7, 9)( 6, 8,11)(12,13,15) \\
                        & b \mapsto  ( 1, 2)( 3, 5, 4, 6)( 7,10,11,14)( 8,12, 9,13)(15,16) \\
n =17 \quad & a \mapsto ( 2, 3, 5)( 6, 8,11)( 7,10, 9)(12,15,13)(14,16,17) \\
       & b \mapsto  ( 1, 2, 4, 7)( 3, 6, 9,12)( 5, 8,10,13)(11,14,15,16).
\end{align*}
This completes the proof of the claim, the lemma, and thus the theorem. \end{proof}
\renewcommand{\qed}{}

\end{proof}

\newpage

\sh{Appendix}

\bigskip

\begin{table}[ht!]
{ \setlength{\tabcolsep}{0.14cm} \tiny
\begin{tabular}{rrrrrrrrrrrrr}
   &  $A_5$  &  $L_2(7)$  &  $A_6$  &  $L_2(8)$  &  $L_2(11)$  &  $L_2(13)$ & $L_2(17)$ & $A_7$  &  $L_2(19)$  &$L_2(16)$ &  $L_3(3)$ &  $U_3(3)$  \\ 
\hline
$A_5$  &  \textbf{1.00}  &  0.02  &  0.13  &  0.05  &  0.17  &  0.03  &  -0.03  &  0.12  &  0.15  &  0.09  &  0.02  &  0.02  \\ 
$L_2(7)$    &  0.02  &  \textbf{1.00}  &  0.04  &  0.23  &  0.05  &  0.16  &  0.05  &  0.06  &  -0.02  &  -0.04  &  0.12  &  0.09  \\ 
$A_6$   &  0.13  &  0.04  &  \textbf{1.00}  &  -0.04  &  0.13  &  -0.07  &  0.02  &  0.10  &  0.11  &  0.09  &  0.04  &  0.00  \\ 
\hline
$L_2(8)$    &  0.05  &  0.23  &  -0.04  &  \textbf{1.00}  &  0.02  &  0.20  &  0.06  &  0.08  &  0.05  &  -0.00  &  -0.00  &  0.11  \\ 
$L_2(11)$   &  0.17  &  0.05  &  0.13  &  0.02  &  \textbf{1.00}  &  -0.01  &  0.03  &  0.11  &  0.11  &  0.14  &  0.07  &  0.05  \\ 
$L_2(13)$   &  0.03  &  0.16  &  -0.07  &  0.20  &  -0.01  &  \textbf{1.00}  &  0.00  &  -0.01  &  0.04  &  0.04  &  0.06  &  0.09  \\ 
\hline
$L_2(17)$  &  -0.03  &  0.05  &  0.02  &  0.06  &  0.03  &  0.00  &  \textbf{1.00}  &  0.01  &  0.05  &  0.03  &  0.11  &  0.12  \\ 
$A_7$  &  0.12  &  0.06  &  0.10  &  0.08  &  0.11  &  -0.01  &  0.01  &  \textbf{1.00}  &  0.08  &  0.10  &  0.03  &  0.11  \\ 
$L_2(19)$  &  0.15  &  -0.02  &  0.11  &  0.05  &  0.11  &  0.04  &  0.05  &  0.08  &  \textbf{1.00}  &  0.11  &  0.03  &  0.03  \\ 
\hline
$L_2(16)$  &  0.09  &  -0.04  &  0.09  &  -0.00  &  0.14  &  0.04  &  0.03  &  0.10  &  0.11  &  \textbf{1.00}  &  -0.02  &  0.07  \\ 
$L_3(3)$  &  0.02  &  0.12  &  0.04  &  -0.00  &  0.07  &  0.06  &  0.11  &  0.03  &  0.03  &  -0.02  &  \textbf{1.00}  &  0.10  \\ 
$U_3(3)$  &  0.02  &  0.09  &  0.00  &  0.11  &  0.05  &  0.09  &  0.12  &  0.11  &  0.03  &  0.07  &  0.10  &  \textbf{1.00}  \\ 
\hline
$L_2(23)$  &  0.01  &  0.10  &  0.03  &  0.07  &  0.05  &  0.03  &  0.12  &  -0.04  &  0.03  &  0.03  &  0.15  &  0.04  \\ 
$L_2(25)$  &  0.04  &  0.06  &  0.15  &  0.06  &  0.14  &  0.03  &  0.13  &  0.09  &  0.10  &  0.10  &  0.21  &  0.08  \\ 
$M_{11}$  &  0.16  &  0.03  &  0.21  &  -0.00  &  0.09  &  -0.02  &  0.09  &  0.12  &  0.01  &  0.05  &  0.05  &  0.06  \\ 
\hline
$L_2(27)$ &  -0.01  &  0.19  &  -0.05  &  0.29  &  0.02  &  0.15  &  0.04  &  0.09  &  0.04  &  0.00  &  0.06  &  0.10  \\ 
$L_2(29)$  &  0.01  &  0.13  &  0.01  &  0.14  &  0.17  &  0.10  &  -0.00  &  0.19  &  0.15  &  0.06  &  0.00  &  0.01  \\ 
$L_2(31)$  &  0.08  &  0.08  &  0.18  &  0.00  &  0.10  &  -0.05  &  0.11  &  0.04  &  0.10  &  0.09  &  0.09  &  0.06  \\ 
\hline
$A_8$  &  0.11  &  0.14  &  0.12  &  0.11  &  0.08  &  0.08  &  0.07  &  0.17  &  0.10  &  0.07  &  0.04  &  0.11  \\ 
$L_3(4)$  &  0.15  &  0.03  &  0.13  &  0.02  &  0.11  &  -0.04  &  0.03  &  0.23  &  0.05  &  0.01  &  0.07  &  0.03  \\ 
$L_2(37)$  &  0.02  &  0.01  &  0.06  &  0.02  &  0.06  &  0.02  &  0.07  &  0.04  &  0.08  &  0.13  &  0.00  &  0.02  \\ 
\hline
$U_4(2)$   &  0.18  &  0.02  &  0.24  &  -0.00  &  0.07  &  -0.04  &  -0.01  &  0.13  &  0.05  &  0.05  &  0.02  &  -0.01  \\ 
Sz(8)   &  -0.00  &  0.02  &  0.11  &  -0.01  &  0.03  &  -0.03  &  0.00  &  -0.02  &  0.09  &  -0.03  &  -0.01  &  -0.03  \\ 
$L_2(32)$   &  0.07  &  0.06  &  -0.02  &  -0.02  &  0.01  &  0.03  &  0.00  &  -0.02  &  0.01  &  0.02  &  -0.00  &  0.05  \\ 
\\ 
\end{tabular}

\vspace{0.6cm}

\begin{tabular}{rrrrrrrrrrrrr}
    &  $L_2(23)$ & $L_2(25)$ & $M_{11}$  &  $L_2(27)$&  $L_2(29)$ & $L_2(31)$&  $A_8$  &  $L_3(4)$  &  $L_2(37)$ & $U_4(2)$  &  Sz(8)  &  $L_2(32)$  \\ 
\hline
$A_5$  &  0.01  &  0.04  &  0.16  &  -0.01  &  0.01  &  0.08  &  0.11  &  0.15  &  0.02  &  0.18  &  -0.00  &  0.07  \\ 
$L_2(7)$  &  0.10  &  0.06  &  0.03  &  0.19  &  0.13  &  0.08  &  0.14  &  0.03  &  0.01  &  0.02  &  0.02  &  0.06  \\ 
$A_6$  &  0.03  &  0.15  &  0.21  &  -0.05  &  0.01  &  0.18  &  0.12  &  0.13  &  0.06  &  0.24  &  0.11  &  -0.02  \\ 
\hline
$L_2(8)$  &  0.07  &  0.06  &  -0.00  &  0.29  &  0.14  &  0.00  &  0.11  &  0.02  &  0.02  &  -0.00  &  -0.01  &  -0.02  \\ 
$L_2(11)$   &  0.05  &  0.14  &  0.09  &  0.02  &  0.17  &  0.10  &  0.08  &  0.11  &  0.06  &  0.07  &  0.03  &  0.01  \\ 
$L_2(13)$ &  0.03  &  0.03  &  -0.02  &  0.15  &  0.10  &  -0.05  &  0.08  &  -0.04  &  0.02  &  -0.04  &  -0.03  &  0.03  \\ 
\hline
$L_2(17)$   &  0.12  &  0.13  &  0.09  &  0.04  &  -0.00  &  0.11  &  0.07  &  0.03  &  0.07  &  -0.01  &  0.00  &  0.00  \\ 
$A_7$  &  -0.04  &  0.09  &  0.12  &  0.09  &  0.19  &  0.04  &  0.17  &  0.23  &  0.04  &  0.13  &  -0.02  &  -0.02  \\ 
$L_2(19)$  &  0.03  &  0.10  &  0.01  &  0.04  &  0.15  &  0.10  &  0.10  &  0.05  &  0.08  &  0.05  &  0.09  &  0.01  \\ 
\hline
$L_2(16)$ &  0.03  &  0.10  &  0.05  &  0.00  &  0.06  &  0.09  &  0.07  &  0.01  &  0.13  &  0.05  &  -0.03  &  0.02  \\ 
$L_3(3)$  &  0.15  &  0.21  &  0.05  &  0.06  &  0.00  &  0.09  &  0.04  &  0.07  &  0.00  &  0.02  &  -0.01  &  -0.00  \\ 
$U_3(3)$  &  0.04  &  0.08  &  0.06  &  0.10  &  0.01  &  0.06  &  0.11  &  0.03  &  0.02  &  -0.01  &  -0.03  &  0.05  \\ 
\hline
$L_2(23)$   &  \textbf{1.00}  &  0.09  &  0.04  &  0.07  &  0.02  &  0.08  &  -0.02  &  0.01  &  0.00  &  0.01  &  -0.04  &  0.08  \\ 
$L_2(25)$   &  0.09  &  \textbf{1.00}  &  0.05  &  0.15  &  0.07  &  0.14  &  0.12  &  0.05  &  0.06  &  0.10  &  0.03  &  0.03  \\ 
$M_{11}$  &  0.04  &  0.05  &  \textbf{1.00}  &  -0.01  &  -0.00  &  0.14  &  0.14  &  0.19  &  0.00  &  0.21  &  0.09  &  0.04  \\ 
\hline
$L_2(27)$   &  0.07  &  0.15  &  -0.01  &  \textbf{1.00}  &  0.19  &  0.01  &  0.11  &  0.02  &  0.03  &  -0.01  &  -0.04  &  0.05  \\ 
$L_2(29)$    &  0.02  &  0.07  &  -0.00  &  0.19  &  \textbf{1.00}  &  0.07  &  0.12  &  0.11  &  0.08  &  0.03  &  -0.01  &  -0.02  \\ 
$L_2(31)$   &  0.08  &  0.14  &  0.14  &  0.01  &  0.07  &  \textbf{1.00}  &  0.09  &  0.10  &  0.02  &  0.13  &  0.08  &  0.08  \\ 
\hline
$A_8$  &  -0.02  &  0.12  &  0.14  &  0.11  &  0.12  &  0.09  &  \textbf{1.00}  &  0.15  &  -0.01  &  0.14  &  0.08  &  -0.03  \\ 
$L_3(4)$   &  0.01  &  0.05  &  0.19  &  0.02  &  0.11  &  0.10  &  0.15  &  \textbf{1.00}  &  -0.00  &  0.21  &  0.26  &  -0.04  \\ 
$L_2(37)$  &  0.00  &  0.06  &  0.00  &  0.03  &  0.08  &  0.02  &  -0.01  &  -0.00  &  \textbf{1.00}  &  0.01  &  -0.03  &  0.06  \\ 
\hline
$U_4(2)$  &  0.01  &  0.10  &  0.21  &  -0.01  &  0.03  &  0.13  &  0.14  &  0.21  &  0.01  &  \textbf{1.00}  &  0.02  &  0.04  \\ 
Sz(8)   &  -0.04  &  0.03  &  0.09  &  -0.04  &  -0.01  &  0.08  &  0.08  &  0.26  &  -0.03  &  0.02  &  \textbf{1.00}  &  -0.03  \\ 
$L_2(32)$   &  0.08  &  0.03  &  0.04  &  0.05  &  -0.02  &  0.08  &  -0.03  &  -0.04  &  0.06  &  0.04  &  -0.03  &  \textbf{1.00}  \\ 
\end{tabular}
}
\vspace{0.2in}
\caption{This table gives the correlations between: (having a cover with group 1, having a cover with group 2).  The average off-diagonal correlation is 0.06.}

\end{table}

\begin{table}[ht!]
{ \setlength{\tabcolsep}{0.14cm} \tiny
\begin{tabular}{rrrrrrrrrrrrr}
  &  $A_5$  &  $L_2(7)$  &  $A_6$  &  $L_2(8)$  &  $L_2(11)$ & $L_2(13)$ & $L_2(17)$&  $A_7$  &  $L_2(19)$&  $L_2(16)$&  $L_3(3)$  &  $U_3(3)$  \\ 
\hline
$A_5$  &  \textbf{1.00}  &  -0.01  &  0.28  &  0.05  &  0.25  &  0.01  &  0.06  &  0.11  &  0.23  &  0.11  &  0.02  &  0.02  \\ 
$L_2(7)$   &  -0.01  &  \textbf{1.00}  &  0.05  &  0.38  &  0.04  &  0.25  &  0.14  &  0.11  &  0.02  &  0.01  &  0.17  &  0.13  \\ 
$A_6$  &  0.28  &  0.05  &  \textbf{1.00}  &  0.00  &  0.22  &  -0.07  &  0.12  &  0.13  &  0.17  &  0.10  &  0.08  &  0.02  \\ 
\hline
$L_2(8)$   &  0.05  &  0.38  &  0.00  &  \textbf{1.00}  &  0.05  &  0.36  &  0.11  &  0.12  &  0.06  &  0.06  &  0.03  &  0.12  \\ 
$L_2(11)$   &  0.25  &  0.04  &  0.22  &  0.05  &  \textbf{1.00}  &  0.03  &  0.07  &  0.06  &  0.18  &  0.12  &  0.08  &  0.04  \\ 
$L_2(13)$  &  0.01  &  0.25  &  -0.07  &  0.36  &  0.03  &  \textbf{1.00}  &  0.07  &  0.01  &  0.04  &  0.10  &  0.08  &  0.13  \\ 
\hline
$L_2(17)$   &  0.06  &  0.14  &  0.12  &  0.11  &  0.07  &  0.07  &  \textbf{1.00}  &  0.07  &  0.12  &  0.07  &  0.15  &  0.11  \\ 
$A_7$    &  0.11  &  0.11  &  0.13  &  0.12  &  0.06  &  0.01  &  0.07  &  \textbf{1.00}  &  0.07  &  0.09  &  0.07  &  0.13  \\ 
$L_2(19)$    &  0.23  &  0.02  &  0.17  &  0.06  &  0.18  &  0.04  &  0.12  &  0.07  &  \textbf{1.00}  &  0.09  &  0.08  &  0.05  \\ \hline
$L_2(16)$   &  0.11  &  0.01  &  0.10  &  0.06  &  0.12  &  0.10  &  0.07  &  0.09  &  0.09  &  \textbf{1.00}  &  0.03  &  0.10  \\ 
$L_3(3)$   &  0.02  &  0.17  &  0.08  &  0.03  &  0.08  &  0.08  &  0.15  &  0.07  &  0.08  &  0.03  &  \textbf{1.00}  &  0.14  \\ 
$U_3(3)$    &  0.02  &  0.13  &  0.02  &  0.12  &  0.04  &  0.13  &  0.11  &  0.13  &  0.05  &  0.10  &  0.14  &  \textbf{1.00}  \\ 
\hline
$L_2(23)$    &  0.06  &  0.13  &  0.02  &  0.04  &  0.06  &  0.05  &  0.13  &  -0.01  &  0.06  &  0.05  &  0.15  &  0.09  \\ 
$L_2(25)$    &  0.12  &  0.13  &  0.20  &  0.14  &  0.17  &  0.06  &  0.17  &  0.12  &  0.14  &  0.15  &  0.21  &  0.13  \\ 
$M_{11}$    &  0.19  &  0.04  &  0.33  &  0.03  &  0.12  &  0.00  &  0.11  &  0.17  &  0.07  &  0.08  &  0.06  &  0.07  \\ 
\hline
$L_2(27)$   &  -0.03  &  0.38  &  -0.06  &  0.45  &  0.05  &  0.35  &  0.06  &  0.10  &  0.01  &  0.01  &  0.09  &  0.16  \\ 
$L_2(29)$   &  0.08  &  0.17  &  0.04  &  0.24  &  0.24  &  0.18  &  0.02  &  0.22  &  0.15  &  0.05  &  0.06  &  0.03  \\ 
$L_2(31)$   &  0.22  &  0.08  &  0.30  &  0.02  &  0.15  &  0.02  &  0.24  &  0.08  &  0.15  &  0.09  &  0.15  &  0.08  \\ 
\hline
$A_8$    &  0.11  &  0.15  &  0.15  &  0.14  &  0.08  &  0.12  &  0.14  &  0.20  &  0.08  &  0.09  &  0.09  &  0.12  \\ 
$L_3(4)$   &  0.21  &  0.08  &  0.27  &  0.04  &  0.15  &  -0.01  &  0.05  &  0.28  &  0.13  &  0.09  &  0.11  &  0.04  \\ 
$L_2(37)$    &  0.09  &  0.03  &  0.12  &  0.05  &  0.14  &  0.10  &  0.15  &  0.02  &  0.09  &  0.16  &  0.03  &  0.08  \\ 
\hline
$U_4(2)$   &  0.17  &  0.03  &  0.34  &  0.01  &  0.10  &  -0.01  &  0.05  &  0.15  &  0.10  &  0.08  &  0.05  &  0.02  \\ 
Sz(8)  &  0.08  &  0.08  &  0.17  &  0.03  &  0.06  &  0.01  &  0.05  &  0.04  &  0.10  &  0.03  &  0.02  &  -0.01  \\ 
$L_2(32)$   &  0.06  &  0.05  &  -0.01  &  0.01  &  0.05  &  0.06  &  0.02  &  -0.01  &  0.04  &  0.05  &  0.01  &  0.08  \\ 
\end{tabular}

\vspace{0.6cm}

\begin{tabular}{rrrrrrrrrrrrr}
   &  $L_2(23)$ &  $L_2(25)$   & $M_{11}$  &  $L_2(27)$  &  $L_2(29)$  &  $L_2(31)$&    $A_8$  &  $L_3(4)$  &  $L_2(37)$ &  $U_4(2)$  &  Sz(8)  &  $L_2(32)$  \\ 
\hline
$A_5$   &  0.06  &  0.12  &  0.19  &  -0.03  &  0.08  &  0.22  &  0.11  &  0.21  &  0.09  &  0.17  &  0.08  &  0.06  \\ 
$L_2(7)$  &  0.13  &  0.13  &  0.04  &  0.38  &  0.17  &  0.08  &  0.15  &  0.08  &  0.03  &  0.03  &  0.08  &  0.05  \\ 
$A_6$  &  0.02  &  0.20  &  0.33  &  -0.06  &  0.04  &  0.30  &  0.15  &  0.27  &  0.12  &  0.34  &  0.17  &  -0.01  \\ 
\hline
$L_2(8)$   &  0.04  &  0.14  &  0.03  &  0.45  &  0.24  &  0.02  &  0.14  &  0.04  &  0.05  &  0.01  &  0.03  &  0.01  \\ 
$L_2(11)$   &  0.06  &  0.17  &  0.12  &  0.05  &  0.24  &  0.15  &  0.08  &  0.15  &  0.14  &  0.10  &  0.06  &  0.05  \\ 
$L_2(13)$   &  0.05  &  0.06  &  0.00  &  0.35  &  0.18  &  0.02  &  0.12  &  -0.01  &  0.10  &  -0.01  &  0.01  &  0.06  \\ 
\hline
$L_2(17)$   &  0.13  &  0.17  &  0.11  &  0.06  &  0.02  &  0.24  &  0.14  &  0.05  &  0.15  &  0.05  &  0.05  &  0.02  \\ 
$A_7$   &  -0.01  &  0.12  &  0.17  &  0.10  &  0.22  &  0.08  &  0.20  &  0.28  &  0.02  &  0.15  &  0.04  &  -0.01  \\ 
$L_2(19)$   &  0.06  &  0.14  &  0.07  &  0.01  &  0.15  &  0.15  &  0.08  &  0.13  &  0.09  &  0.10  &  0.10  &  0.04  \\ 
\hline
$L_2(16)$    &  0.05  &  0.15  &  0.08  &  0.01  &  0.05  &  0.09  &  0.09  &  0.09  &  0.16  &  0.08  &  0.03  &  0.05  \\ 
$L_3(3)$    &  0.15  &  0.21  &  0.06  &  0.09  &  0.06  &  0.15  &  0.09  &  0.11  &  0.03  &  0.05  &  0.02  &  0.01  \\ 
$U_3(3)$    &  0.09  &  0.13  &  0.07  &  0.16  &  0.03  &  0.08  &  0.12  &  0.04  &  0.08  &  0.02  &  -0.01  &  0.08  \\ 
\hline
$L_2(23)$    &  \textbf{1.00}  &  0.11  &  0.05  &  0.04  &  0.04  &  0.08  &  0.01  &  0.04  &  0.09  &  0.02  &  -0.05  &  0.15  \\ 
$L_2(25)$    &  0.11  &  \textbf{1.00}  &  0.12  &  0.15  &  0.15  &  0.18  &  0.20  &  0.06  &  0.16  &  0.14  &  0.04  &  0.02  \\ 
$M_{11}$   &  0.05  &  0.12  &  \textbf{1.00}  &  -0.04  &  0.02  &  0.22  &  0.14  &  0.24  &  0.05  &  0.25  &  0.08  &  0.05  \\ 
\hline
$L_2(27)$   &  0.04  &  0.15  &  -0.04  &  \textbf{1.00}  &  0.25  &  -0.03  &  0.10  &  0.02  &  0.03  &  0.01  &  0.01  &  0.08  \\ 
$L_2(29)$  &  0.04  &  0.15  &  0.02  &  0.25  &  \textbf{1.00}  &  0.06  &  0.13  &  0.21  &  0.12  &  0.04  &  0.04  &  0.01  \\ 
$L_2(31)$  &  0.08  &  0.18  &  0.22  &  -0.03  &  0.06  &  \textbf{1.00}  &  0.10  &  0.12  &  0.12  &  0.15  &  0.10  &  0.06  \\ 
\hline
$A_8$   &  0.01  &  0.20  &  0.14  &  0.10  &  0.13  &  0.10  &  \textbf{1.00}  &  0.18  &  0.07  &  0.16  &  0.09  &  -0.03  \\ 
$L_3(4)$  &  0.04  &  0.06  &  0.24  &  0.02  &  0.21  &  0.12  &  0.18  &  \textbf{1.00}  &  0.02  &  0.25  &  0.30  &  -0.04  \\ 
$L_2(37)$    &  0.09  &  0.16  &  0.05  &  0.03  &  0.12  &  0.12  &  0.07  &  0.02  &  \textbf{1.00}  &  0.06  &  0.03  &  0.10  \\ \hline
$U_4(2)$    &  0.02  &  0.14  &  0.25  &  0.01  &  0.04  &  0.15  &  0.16  &  0.25  &  0.06  &  \textbf{1.00}  &  -0.01  &  0.01  \\ 
Sz(8)    &  -0.05  &  0.04  &  0.08  &  0.01  &  0.04  &  0.10  &  0.09  &  0.30  &  0.03  &  -0.01  &  \textbf{1.00}  &  -0.07  \\ 
$L_2(32)$    &  0.15  &  0.02  &  0.05  &  0.08  &  0.01  &  0.06  &  -0.03  &  -0.04  &  0.10  &  0.01  &  -0.07  &  \textbf{1.00}  \\ 
\end{tabular}

}
\vspace{0.2in}
\caption{This table gives the correlations between: (having a cover with group 1 with positive betti number, having a cover with group 2 with positive betti number).  The average off-diagonal correlation is 0.09.}
\end{table}

\end{document}

%% file: gtoutput.tex
\def\ifplaintex{\expandafter\ifx\csname documentclass\endcsname\relax}

\ifplaintex 
\else
\fi

\expandafter\ifx\csname beginpicture\endcsname\relax
\expandafter\ifx\csname documentclass\endcsname\relax
\input pictex \else
\input prepictex \input pictex \input postpictex \fi\fi

\def\gt{{\mathsurround=0pt\it $\cal G\mskip-2mu$eometry \&\ 
$\cal T\!\!$opology}}        

\def\gtp{{\mathsurround=0pt\it $\cal G\mskip-2mu$eometry \&\ 
$\cal T\!\!$opology $\cal P\!$ublications}}  

\def\lognumber#1{\def\thelognumber{#1}}
\def\volumenumber#1{\def\thevolumenumber{#1}}
\def\papernumber#1{\def\thepapernumber{#1}}
\def\volumeyear#1{\def\thevolumeyear{#1}}

\def\pagenumbers#1#2{\def\startpage{#1}\def\finishpage{#2}}
\def\published#1{\def\publishdate{#1}}
\def\proposed#1{\def\theproposer{#1}}
\def\seconded#1{\def\theseconders{#1}}
\def\received#1{\def\receiveddate{#1}}

\def\accepted#1{\def\accepteddate{#1}}

\def\asciiaddress#1{\def\theasciiaddress{#1}}
\def\asciiemail#1{\def\theasciiemail{#1}}

\def\shortauthors#1{\def\theshortauthors{#1}}

\let\\\par\let\thelognumber\relax
\let\thevolumenumber\relax\let\thepapernumber\relax
\let\thevolumeyear\relax\let\thesamplenumber\relax\let\startpage\relax
\let\finishpage\relax\let\publishdate\relax\let\receiveddate\relax
\let\reviseddate\relax\let\accepteddate\relax\let\theasciititle\relax
\let\theasciiauthors\relax\let\theasciiaddress\relax
\let\theasciiabstract\relax
\let\theasciiemail\relax\let\theshortauthors\relax\let\theshorttitle\relax

\long\def\maketitlep{   

\count0=\startpage

\gt\hfill      
\beginpicture
\setcoordinatesystem units <0.33truein, 0.33truein> point at 2.2 0.9
\setplotsymbol ({$\cal G$})
\plotsymbolspacing=9truept
\circulararc 315 degrees from 0 1 center at 0 0
\setplotsymbol ({$\cal T$})
\circulararc 315 degrees from 1 -1 center at 1 0
\endpicture
\break
{\small\ifx\thesamplenumber\relax 
Volume \else Sample
\fi\thevolumenumber\ (\thevolumeyear)
\startpage--\finishpage\nl
Published: \publishdate}
\vglue 0.5truein plus 0.4fil minus 0.1truein

{\parskip=0pt\leftskip 0pt plus 1fil\def\\{\par\smallskip}{\ifplaintex\large
\else\Large\fi\bf\thetitle}\par\medskip}   

\vglue 0pt plus 0.1fil 

{\parskip=0pt\leftskip 0pt plus 1fil\def\\{\par}{\sc\theauthors}
\par\medskip}

\vglue 0pt plus 0.1fil 

{\small\parskip=0pt\let\newline\\
{\leftskip 0pt plus 1fil\def\\{\par}{\sl\theaddress}\par}
\expandafter\ifx\theemail\relax    
\relax\else\vglue 5pt plus 0.02fil minus 2pt\def\\{\stdspace{\rm 
and}\stdspace} 
\cl{Email:\stdspace\tt\theemail}\fi
\ifx\theurl\relax                  
\relax\else\vglue 5pt plus 0.02fil minus 2pt\def\\{\stdspace{\rm 
and}\stdspace}
\cl{URL:\stdspace\tt\theurl}\fi\par}

\vglue 7pt plus 0.3fil minus 3pt

{\bf Abstract}
\vglue 5pt plus 0.1fil minus 2pt

\theabstract

\vglue 7pt plus 0.3fil minus 3pt

{\bf AMS Classification numbers}\quad Primary:\quad \theprimaryclass

Secondary:\quad \thesecondaryclass

\vglue 5pt plus 0.3fil minus 2pt

{\bf Keywords}\quad \thekeywords

\vglue 10pt plus 0.5fil minus 5pt

{\small  Proposed: \theproposer\hfill Received: \receiveddate\nl
Seconded: \theseconders\hfill 
\ifx\reviseddate\relax                         
Accepted: \accepteddate                        
\else
Revised: \reviseddate                          
\fi}
\eject
}       

\let\maketitlepage\maketitlep
\let\maketitle\maketitlepage

\font\phead=cmsl9 scaled 950
\font=cmsl9 scaled 1050
\font\pnum=cmbx10 scaled 913
\font\lnum=cmbx10 
\font\pfoot=cmsl9 scaled 950
\font=cmsl9 scaled 1050
\ifplaintex
\headline{\vbox to 0pt{\vskip -4.5mm\line{\small\phead\ifnum
\count0=\startpage ISSN 1364-0380 (on line)
1465-3060 (printed) \hfill {\pnum\folio}\else\ifodd\count0\def\\{ }%
\ifx\theshorttitle\relax\thetitle\else\theshorttitle\fi\hfill{\pnum\folio}
\else\def\\{ and }{\pnum\folio}\hfill\ifx\theshortauthors\relax\theauthors
\else\theshortauthors\fi\fi\fi}\vss}}
\footline{\vbox to 0pt{\vglue 0mm\line{\small\pfoot\ifnum\count0=\startpage
\copyright\ \gtp\hfill\else
\gt, Volume \thevolumenumber\ (\thevolumeyear)\hfill\fi}\vss
}}
\else
\makeatletter
\def\@oddhead{{\small\lhead\ifnum\count0=\startpage ISSN 1364-0380 (on line)
1465-3060 (printed) \hfill {\lnum\number\count0}\else\ifodd\count0
\def\\{ }\ifx\theshorttitle\relax \thetitle \else\theshorttitle\fi\hfill
{\lnum\number\count0}\else\def\\{ and }{\lnum\number\count0}
\hfill\ifx\theshortauthors\relax 
\theauthors\else\theshortauthors\fi\fi\fi}}\def\@evenhead{\@oddhead}
\def\@oddfoot{\small\lfoot\ifnum\count0=\startpage\copyright\ \gtp\hfill\else
\gt, Volume \thevolumenumber\ (\thevolumeyear)\hfill\fi}
\def\@evenfoot{\@oddfoot}
\makeatother
\fi

\newwrite\gtoutfile
\long\gdef\makeheadfile{  
{\def\\{, }\def\s{ }
\immediate\openout\gtoutfile head.xxx
\immediate\write\gtoutfile{To: math@arxiv.org}
\immediate\write\gtoutfile{Subject: put or rep NNNNN:pppp}
\immediate\write\gtoutfile{--text follows this line--}
\immediate\write\gtoutfile{Proxy-for: \ifx\theasciiauthors\relax
\theauthors\else\theasciiauthors\fi\s<\ifx\theasciiemail\relax\theemail\else\theasciiemail\fi>}
\immediate\write\gtoutfile{\noexpand\\}
\immediate\write\gtoutfile{Authors: \ifx\theasciiauthors\relax
\theauthors\else\theasciiauthors\fi}
{\def\\{ }\immediate\write\gtoutfile{Title: \ifx\theasciititle\relax
\thetitle\else\theasciititle\fi}}
\immediate\write\gtoutfile{Subj-class: GT or SG or MG etc}
\immediate\write\gtoutfile{MSC-class: \theprimaryclass\ifx\thesecondaryclass\relax\else, \thesecondaryclass\fi}
\immediate\write\gtoutfile{Journal-ref: Geom. Topol. \thevolumenumber
(\thevolumeyear) \startpage-\finishpage}
\immediate\write\gtoutfile{Comments: Published by Geometry and Topology at}
\immediate\write\gtoutfile{\s\s http://www.maths.warwick.ac.uk/gt/GTVol\thevolumenumber/paper\thepapernumber.abs.html}
\immediate\write\gtoutfile{\noexpand\\}
\immediate\write\gtoutfile{}
\ifx\theasciiabstract\relax
\immediate\write\gtoutfile{\theabstract}\else
\immediate\write\gtoutfile{\theasciiabstract}\fi
\immediate\write\gtoutfile{}
\immediate\write\gtoutfile{\noexpand\\}
\immediate\write\gtoutfile{}
\immediate\closeout\gtoutfile}}  

\def\maketitlepage{\maketitlep\makeheadfile}
\let\maketitle\maketitlepage